\documentclass[11pt,reqno]{amsart}
\usepackage[pdfstartview=FitB]{hyperref}
\usepackage{natbib}
\usepackage{graphicx}
\usepackage{geometry}                
\usepackage{amsmath,amssymb,amsthm}

\theoremstyle{plain}

\newtheorem{algorithm}{Algorithm}
\newtheorem{lemma}{Lemma}

\newtheorem{proposition}{Proposition}
\theoremstyle{definition}

\DeclareMathOperator{\Cov}{Cov}
\DeclareMathOperator{\Var}{Var}

\newcommand{\citen}{\citeasnoun}
\newcommand{\Pb}{\mathbb{P}}

\newcommand{\G}{\Gamma}

\renewcommand{\qed}{\hfill{\tiny \ensuremath{\blacksquare} }}%

\newcommand{\Ep}{{\mathrm{E}}}

\newcommand{\En}{{\mathbb{E}_n}}
\renewcommand{\Pr}{{\mathrm{P}}}

\newcommand{\mB}{\mathcal{B}}

\newcommand{\M}{\mathrm{M}}
\newcommand{\A}{\mathrm{A}}

\newcommand{\spt}{\mathrm{support}}
\newcommand{\sC}{\mathsf{C}}
\renewcommand{\sc}{\mathsf{c}}
\newcommand{\sA}{\mathsf{A}}
\newcommand{\sa}{\mathsf{a}}
\newcommand{\sq}{\mathsf{q}}

\renewcommand{\AA}{\mathcal{A}}

\newcommand{\Mn}{\hat \M}

\newcommand{\citeasnoun}{\cite}
\newcommand{\Ln}{\widehat{\mathrm{L}}}

\begin{document}
\begin{center}
\bigskip
{\huge Valid Post-Selection and Post-Regularization Inference: An Elementary, General Approach}\\

\end{center}

\title[Post-Selection and Post-Regularization Inference]{}
\author{Victor Chernozhukov, Christian Hansen, and Martin Spindler}\thanks{Chernozhukov: Massachussets Institute of Technology, 50 Memorial Drive, E52-361B, Cambridge, MA 02142, vchern@mit.edu.
Hansen: University of Chicago Booth School of Business, 5807 S. Woodlawn Ave., Chicago, IL 60637, chansen1@chicagobooth.edu. Spindler:  Munich Center for the Economics of Aging, Amalienstr. 33, 80799 Munich, Germany, spindler@mea.mpisoc.mpg.de. }
\date{December, 2014. We thank Denis Chetverikov,  Mert Demirer, Anna Mikusheva, seminar participants and the discussant Susan Athey at
the AEA Session on Machine Learning in Economics and Econometrics, CEME conference on Non-Standard Problems in Econometrics,
Berlin Statistics Seminar, and the students from MIT's 14.387 Applied Econometrics Class for useful comments.}
\maketitle

\begin{footnotesize}
\textbf{Abstract.} We present an expository, general analysis of valid post-selection or post-regularization inference about a low-dimensional target parameter in the presence of a very high-dimensional nuisance parameter that is estimated using selection or regularization methods. Our analysis provides a set of high-level conditions under which inference for the low-dimensional parameter based on testing or point estimation methods will be regular despite selection or regularization biases occurring in the estimation of the high-dimensional nuisance parameter. The results may be applied to establish uniform validity of post-selection or post-regularization inference procedures for low-dimensional target parameters over large classes of models.  The high-level conditions allow one to clearly see the types of structure needed to achieve valid post-regularization inference and encompass many existing results.  A key element of the structure we employ and discuss in detail is the use of so-called orthogonal or ``immunized'' estimating equations that are locally insensitive to small mistakes in estimation of the high-dimensional nuisance parameter.  As an illustration, we use the high-level conditions to provide readily verifiable sufficient conditions for a class of affine-quadratic models that include the usual linear model and linear instrumental variables model as special cases.  As a further application and illustration, we use these results to provide an analysis of post-selection inference in a  linear instrumental variables model with many regressors and many instruments.  We conclude with a review of other developments in post-selection inference and note that many of the developments can be viewed as special cases of the general encompassing framework of orthogonal estimating equations provided in this paper.


\end{footnotesize}


\section{Introduction}

Analysis of high-dimensional models, models in which the number of parameters to be estimated is large relative to the sample size, is becoming increasingly important.  Such models arise naturally in readily available high-dimensional data which have many measured characteristics available per individual observation as in, for example, large survey data sets, scanner data, and text data.  Such models also arise naturally even in data with a small number of measured characteristics in situations where the exact functional form with which the observed variables enter the model is unknown.  Examples of this scenario include semiparametric models with nonparametric nuisance functions.  More generally, models with many parameters relative to the sample size often arise when attempting to model complex phenomena.  

The key concept underlying the analysis of high-dimensional models is that regularization, such as model selection or shrinkage of model parameters, is necessary if one is to draw meaningful conclusions from the data.  For example, the need for regularization is obvious in a linear regression model with the number of right-hand-side variables greater than the sample size, but arises far more generally in any setting in which the number of parameters is not small relative to the sample size.  Given the importance of the use of regularization in analyzing high-dimensional models, it is then important to explicitly account for the impact of this regularization on the behavior of estimators if one wishes to accurately characterize their finite-sample behavior. The use of such regularization techniques may easily invalidate conventional approaches to inference about model parameters and other interesting target parameters.  A major goal of this paper is to present a general, formal framework that provides guidance about setting up estimating equations and making appropriate use of regularization devices so that inference about parameters of interest will remain valid in the presence of data-dependent model selection or other approaches to regularization. 

It is important to note that understanding estimators' behavior in high-dimensional settings is also useful in conventional low-dimensional settings.  As noted above, dealing formally with high-dimensional models requires that one explicitly accounts for model selection or other forms of regularization.  Providing results that explicitly account for this regularization then allows us to accommodate and coherently account for the fact that low-dimensional models estimated in practice are often the result of specification searches.  As in the high-dimensional setting, failure to account for this variable selection will invalidate the usual inference procedures, whereas the approach that we outline will remain valid and can easily be applied in conventional low-dimensional settings.


The chief goal of this overview paper is to offer a general framework that encompasses many existing results regarding inference on model parameters in high-dimensional models.  The encompassing framework we present and the key theoretical results are new, although they are clearly heavily influenced and foreshadowed by previous, more specialized results.  As an application of the framework, we also present new results on inference in a reasonably broad class of models, termed affine-quadratic models, that includes the usual linear model and linear instrumental variables (IV) model and then apply these results to provide new ones regarding post-regularization inference on the parameters on endogenous variables in a linear instrumental variables model with very many instruments and controls (and also allowing for some misspecification).  We also provide a discussion of previous research that aims to highlight that many existing results fall within the general framework.

Formally, we present a series of results for obtaining valid inferential statements about a low-dimensional parameter of interest, $\alpha$, in the presence of a high-dimensional nuisance parameter $\eta$. The general approach we offer relies on two fundamental elements.  First, it is important that estimating equations used to draw inferences about $\alpha$ satisfy a key orthogonality or immunization condition.\footnote{We refer to the condition as an orthogonality or immunization condition as orthogonality is a much used term and our usage differs from some other usage in defining orthogonality conditions used in econometrics.}  For example, when estimation and inference for $\alpha$ are based on the empirical analog of a theoretical system of equations 
$$\M(\alpha, \eta)=0,$$ 
we show that setting up the equations in a manner such that the orthogonality or immunization condition 
$$\partial_\eta \M(\alpha, \eta) = 0$$
holds is an important element in providing an inferential procedure for $\alpha$ that remains valid when $\eta$ is estimated using regularization.  We note that this condition can generally be established.  For example, we can apply Neyman's classic orthogonalized score in likelihood settings; see, e.g. \citen{Neyman59} and \citen{Neyman1979}.  We also describe an extension of this classic approach to the GMM setting.  In general, applying this orthogonalization will introduce additional nuisance parameters that will be treated as part of $\eta$.   

The second key element of our approach is the use of high-quality, structured estimators of $\eta$.  Crucially, additional structure on $\eta$ is needed for informative inference to proceed, and it is thus important to use estimation strategies that leverage and perform well under the desired structure.  An example of a structure that has been usefully employed in the recent literature is approximate sparsity, e.g. \citen{BCCH12}.  Within this framework, $\eta$ is well approximated by a sparse vector which suggests the use of a sparse estimator such as the Lasso (\citen{FF:1993} and \citen{T1996}).  
The Lasso estimator solves the general problem
$$
\widehat\eta_L = \arg\min_{\eta } \ \ell(\mathrm{data}, \eta) + \lambda\sum_{j=1}^{p} |\psi_j\eta_j|,
$$
where $\ell(\mathrm{data}, \eta)$ is some general loss function that depends on the data and the parameter $\eta$, $\lambda$ is a penalty level, and $\psi_j$'s are penalty loadings.  The leading example is the usual linear model in which $\ell(\mathrm{data},\eta) = \sum_{i=1}^{n} (y_i - x_i'\eta)^2$ is the usual least-squares loss, with $y_i$ denoting the outcome of interest for observation $i$ and $x_i$ denoting predictor variables, and we provide further discussion of this example in the appendix.  Other examples of $\ell(\mathrm{data}, \eta)$ include suitable loss functions corresponding to well-known M-estimators, the negative of the log-likelihood, and GMM criterion functions.  This estimator and related methods such as those in \citen{CandesTao2007}, \citen{MY2007}, \citen{BickelRitovTsybakov2009}, \citen{BC-PostLASSO}, and \citen{BCW-SqLASSO} are computationally efficient and have been shown to have good estimation properties even when perfect variable selection is not feasible under approximate sparsity.  These good estimation properties then translate into providing ``good enough'' estimates of $\eta$ to result in valid inference about $\alpha$ when coupled with orthogonal estimating equations as discussed above.  
Finally, it is important to note that the general results we present do not require or leverage approximate sparsity or sparsity-based estimation strategies.  We provide this discussion here simply as an example and because the structure offers one concrete setting in which the general results we establish may be applied.

In the remainder of this paper, we present the main results.  In Sections 2 and 3, we provide our general set of results that may be used to establish uniform validity of inference about low-dimensional parameters of interest in the presence of high-dimensional nuisance parameters.  We provide the framework in Section 2, and then discuss how to achieve the key orthogonality condition in Section 3.  In Sections 4 and 5, we provide details about establishing the necessary results for the estimation quality of $\eta$ within the approximately sparse framework.  The analysis in Section 4 pertains to a reasonably general class of affine-quadratic models, and the analysis of Section 5 specializes this result to the case of estimating the parameters on a vector of endogenous variables in a linear instrumental variables model with very many potential control variables and very many potential instruments.  The analysis in Section 5 thus extends results from \citen{BCCH12} and \citen{BelloniChernozhukovHansen2011}.  We also provide a brief simulation example and an empirical example that looks at logit demand estimation within the linear many instrument and many control setting in Section 5.  We conclude with a literature review in Section 6.





\textbf{Notation.} We use ``wp $\to 1$" to abbreviate the phrase ``with probability
that converges to 1", and we use the arrows $\to_{\Pr_n}$ and $\leadsto_{\Pr_n}$ to denote convergence in
probability and in distribution under the sequence of probability measures $\{\Pr_n\}$.  The symbol $\sim$
means ``distributed as". The notation $a \lesssim b$ means that $a = O(b)$ and $a \lesssim_{\Pr_n} b$ means that $a = O_{\Pr_n}(b)$.
The $\ell_{2}$ and $\ell_{1}$ norms are denoted by
$\|\cdot\|$ and $\| \cdot \|_{1}$, respectively;  and the $\ell_{0}$-``norm", $\|\cdot\|_0$, denotes the number of non-zero components of a vector.  When applied to a matrix, $\|\cdot\|$ denotes the operator norm.
We use the notation $a \vee b = \max( a, b)$ and $a \wedge b = \min(a , b)$. Here and below, $\En[\cdot]$ abbreviates the average $n^{-1}\sum_{i=1}^n[\cdot]$ over index $i$.  That is, $\En[f(w_i)]$ denotes $n^{-1}\sum_{i=1}^n[f(w_i)]$. In what follows, we use the $m$-sparse norm of a matrix $Q$ defined as
$$
\|Q\|_{\mathsf{sp}(m)} = \sup\{ | b'Qb|/\|b\|^2 :  \|b\|_0 \leq m, \|b\| \neq 0\}.
$$
We also consider the pointwise norm of a square matrix matrix $Q$ at a point $x \neq 0$:
$$
\|Q\|_{\mathsf{pw}(x)} = |x'Qx|/\|x\|^2.
$$
For a differentiable map $x \mapsto f(x)$, mapping $\mathbb{R}^d$ to $\mathbb{R}^k$,  we use $\partial_{x'} f$ to abbreviate the partial derivatives  $(\partial/\partial x') f$, and we correspondingly use the expression
$\partial_{x'} f(x_0)$ to mean $\partial_{x'} f (x) \mid_{x = x_0}$, etc.  We use $x'$ to denote the transpose of a column vector $x$.

\section{A Testing and Estimation Approach to Valid Post-Selection and Post-Regularization Inference}

\subsection{The Setting} We assume that estimation is based on the first $n$ elements $(w_{i,n})_{i=1}^n$ of the \textit{stationary} data-stream $(w_{i,n})_{i=1}^\infty$ which lives on the probability space $(\Omega, \mathcal{A}, \Pr_n)$.   The data points
 $w_{i,n}$ take values in a measurable space $\mathcal{W}$ for each $i$ and $n$. Here, $\Pr_n$, the probability law or data-generating process, can change with $n$.  We allow the law to change with $n$ to claim robustness or uniform validity of results with respect to perturbations of such laws.  Thus the data, all parameters, estimators, and other quantities are indexed by $n$, but we typically suppress this dependence to simplify notation.

The target parameter value $\alpha=\alpha_0$ is assumed to solve the system of theoretical equations
$$
{\M}(\alpha, \eta_0) = 0,  $$
where $\M= (\M_l)_{l=1}^k$ is a measurable map from $\mathcal{A}\times \mathcal{H}$
to $\mathbb{R}^{k}$ and $\mathcal{A}\times \mathcal{H}$ are some convex subsets of $\mathbb{R}^d \times \mathbb{R}^p$. Here the dimension $d$ of the target parameter $\alpha \in \mathcal{A}$ and the number of equations $k$ are assumed to be fixed and
the dimension $p=p_n$ of the nuisance parameter $\eta \in \mathcal{H}$ is allowed to be very high, potentially much larger than $n$.
To handle the high-dimensional nuisance parameter $\eta$, we employ structured assumptions and selection or
regularization methods appropriate for the structure to estimate $\eta_0$.

Given an appropriate estimator $\hat \eta$, we can
construct an estimator $\hat \alpha$ as an approximate solution 
to the estimating equation:
$$
\| \hat{\M}(\hat \alpha, \hat \eta)\| \leq \inf_{\alpha \in \AA} \| \hat{\M}(\alpha, \hat \eta)\| + o(n^{-1/2})
$$
where $\hat{\M} = (\hat{\M}_l)_{l=1}^k$ is the empirical analog of theoretical equations $\M$, which is a measurable map from
$\mathcal{W}^n\times \mathcal{A}\times \mathcal{H}$ to $\mathbb{R}^k$.  
We can also use  $\hat{\M}(\alpha, \hat \eta)$ to test hypotheses about $\alpha_0$ and then
invert the tests to construct confidence sets.

It is not required in the formulation above, but a typical case is when $\hat \M$ and $\M$
are formed as theoretical and empirical moment functions:
$$
 \quad   {\M}(\alpha, \eta) :=  \Ep[\psi(w_i, \alpha, \eta)],
 \quad   \hat{\M}(\alpha, \eta) :=  \En[\psi(w_i, \alpha, \eta)],
 $$
where $\psi = (\psi_l)_{l=1}^k$ is a measurable map from $\mathcal{W}\times \mathcal{A}\times \mathcal{H}$
to $\mathbb{R}^{k}$. Of course, there are many problems that do not fall in the moment condition framework.

\subsection{Valid Inference via Testing}
A simple introduction to the inferential problem is via the testing problem in which
we would like to test some hypothesis about the true parameter value $\alpha_0$.
By inverting the test, 
we create a confidence set for $\alpha_0$. The key condition
for the validity of this confidence region is adaptivity, which can be ensured
by using orthogonal estimating equations and using structured assumptions on the high-dimensional
nuisance parameter.\footnote{We refer to \cite{Bickel1982} for a definition of and introduction to adaptivity.}

The key condition enabling us to perform valid inference on $\alpha_0$
is the \textit{adaptivity} condition:
\begin{equation}\label{eq:adaptivity}
\sqrt{n}(\hat{\M}(\alpha_0, \hat \eta) - \hat{\M}(\alpha_0, \eta_0))    \to_{\Pr_n} 0.
\end{equation}
This condition states that using $\sqrt{n}\hat{\M}(\alpha_0, \hat \eta)$ is as good as
using $\sqrt{n}\hat{\M}(\alpha_0, \eta_0)$, at least to the first order.  This condition may hold
despite using estimators $\hat \eta$ that are not asymptotically
linear and are  non-regular.    Verification of adaptivity may involve substantial work as illustrated below.
A key requirement that often arises is the \textit{orthogonality} or \textit{immunization} condition:
\begin{equation}\label{eq:ortho}
\partial_{\eta'} {\M}(\alpha_0, \eta_0) = 0.
\end{equation}
This condition states that the equations are locally insensitive to small perturbations
of the nuisance parameter around the true parameter values.  In several important models, this condition
is equivalent to the double-robustness condition 
(\cite{robins:dr}). Additional assumptions regarding the \textit{quality of
estimation} of $\eta_0$ are also needed and are highlighted below.

The adaptivity condition immediately allows
us to use the  statistic $\sqrt{n}\hat{\M}(\alpha_0, \hat \eta)$ to perform inference.
Indeed, suppose we have that
\begin{equation}\label{eq:normality}
 \Omega^{-1/2}(\alpha_0)\sqrt{n}\hat{\M}(\alpha_0, \eta_0) \leadsto_{\Pr_n} \mathcal{N}(0,I_k)
\end{equation}
for some positive definite $\Omega(\alpha)  = \Var( \sqrt{n}\hat{\M}(\alpha, \eta_0))$.
This condition can be verified  using central limit theorems for triangular arrays. Such
theorems are available for independently and identically distributed (i.i.d.) as well as dependent and clustered data. Suppose further that there exists $\hat \Omega(\alpha)$ such that
\begin{equation}\label{eq:variance}
\hat \Omega^{-1/2}(\alpha_0) \Omega^{1/2}(\alpha_0) \to_{\Pr_n} I_k.
\end{equation}
It is then immediate that the following score statistic, evaluated at $\alpha = \alpha_0$, is asymptotically normal,
\begin{equation}\label{eq:normality2}
S(\alpha) :=  \hat \Omega^{-1/2}_n(\alpha)\sqrt{n}\hat{\M}(\alpha, \hat \eta) \leadsto_{\Pr_n} \mathcal{N}(0,I_k),
\end{equation}
and that the quadratic form of this score statistic is asymptotically $\chi^2$
with $k$ degrees of freedom:
\begin{equation}\label{eq:chi-square}
C(\alpha_0) = \| S(\alpha_0)\|^2 \leadsto_{\Pr_n}  \chi^2(k).
\end{equation}

The statistic given in (\ref{eq:chi-square}) simply corresponds to a quadratic form in appropriately normalized statistics that have the desired immunization or orthogonality condition.  We refer to this statistic as a ``generalized $C(\alpha)$-statistic'' in honor of Neyman's  fundamental contributions, e.g. \citen{Neyman59} and \citen{Neyman1979}, because, in likelihood settings, the statistic (\ref{eq:chi-square}) reduces to Neyman's $C(\alpha)$-statistic and the generalized score $S(\alpha_0)$ given in (\ref{eq:normality2}) reduces to Neyman's orthogonalized score.  We demonstrate these relationships in the special case of likelihood models in Section 3.1 and provide a generalization to GMM models in Section 3.2. Both of these examples serve to illustrate the construction of appropriate statistics in different settings, but we note that the framework applies far more generally.


The following elementary result is an immediate consequence of the preceding discussion. \\

\begin{proposition}[Valid Inference After Selection or Regularizaton] \textit{Consider a sequence $\{\mathbf{P}_n\}$ of sets
of probability laws such that for each sequence $\{\Pr_n\} \in \{\mathbf{P}_n\}$ the
adaptivity condition (\ref{eq:adaptivity}), the normality condition (\ref{eq:normality}), and the variance consistency condition (\ref{eq:variance}) hold.
Then $\mathsf{CR}_{1-a} = \{ \alpha \in \mathcal{A}:  C(\alpha) \leq  c(1-a) \}$, where
$c(1-a)$ is the $1-a$-quantile of a $\chi^2(k)$, is a uniformly valid confidence
interval for $\alpha_0$ in the sense that}
$$
\lim_{n \to \infty} \sup_{\Pr \in \mathbf{P}_n}  |\Pr (  \alpha_0 \in \mathsf{CR}_{1-a} ) - (1-a) |  = 0.
$$
\end{proposition}

We remark here that in order to make the uniformity claim interesting we should insist that
the sets of probability laws $\mathbf{P}_n$ are non-decreasing in $n$, i.e.
$\mathbf{P}_{\bar n} \subseteq \mathbf{P}_{n}$ whenever $\bar n \leq n$.

Proof. For any sequence of positive constants $\epsilon_n$ approaching $0$, let
 $\Pr_n \in \mathbf{P}_n$ be any sequence such that
 $$
  |\Pr_n (  \alpha_0 \in \mathsf{CR}_{1-a} ) - (1-a) | +\epsilon_n \geq \sup_{\Pr \in \mathbf{P}_n}  |\Pr (  \alpha_0 \in \mathsf{CR}_{1-a} ) - (1-a) |.
$$
By conditions (\ref{eq:normality}) and (\ref{eq:variance}) we have that
$$
\Pr_n (  \alpha_0 \in \mathsf{CR}_{1-a} )  =\Pr_n (  C(\alpha_0) \leq c(1-a) ) \to \mathbb{P} (\chi^2(k) \leq c(1-a)) = 1-a,
$$
which implies the conclusion from the preceding display. \qed

\subsection{Valid Inference via Adaptive Estimation}   Suppose that $\M(\alpha_0, \eta_0) =0$ holds 
for $\alpha_0 \in \mathcal{A}$.  We consider an estimator $\hat \alpha \in \mathcal{A}$
that is an approximate minimizer of the map $\alpha \mapsto \|  \hat{\M}(\alpha, \hat \eta)\|$ in the sense that
\begin{equation}\label{eq:minimization}
\| \hat \M(\hat \alpha, \hat \eta)\| \leq \inf_{\alpha \in \mathcal{A}} \| \hat \M(\alpha, \hat \eta)\| + o(n^{-1/2}).
\end{equation}

In order to analyze this estimator, we
assume that the derivatives
$\Gamma_1 : =  \partial_{\alpha'} \M(\alpha_0, \eta_0)$ and $\partial_{\eta'} \M(\alpha, \eta_0)$
exist.  We assume that $\alpha_0$ is interior relative
to the parameter space $\mathcal{A}$; namely, for some $\ell_n \to \infty$ such that $\ell_n/\sqrt{n} \to 0$,
\begin{equation}\label{eq:interior}
\{ \alpha \in \mathbb{R}^d: \| \alpha - \alpha_0\| \leq \ell_n/\sqrt{n} \} \subset \mathcal{A}.
\end{equation}
We also assume that the following local-global identifiability condition holds: For some constant $c>0$,
\begin{equation}
2 \| \M(\alpha, \eta_0)\| \geq  \| \Gamma_1 (\alpha - \alpha_0) \| \wedge c  \ \ \forall \alpha \in \mathcal{A}, \quad \mathrm{mineig}{(\Gamma_1'\Gamma_1)}  \geq c.\label{eq:id}
\end{equation}
Furthermore, for $\Omega = \Var ( \sqrt{n} \hat \M(\alpha_0,  \eta_0) )$, we suppose that  the 
central limit theorem, 
\begin{equation}\label{eq:normal score}
\Omega^{-1/2} \sqrt{n} \hat \M(\alpha_0, \eta_0) \leadsto_{\Pr_n} \mathcal{N}(0, I),
 \end{equation}
and the stability condition,
\begin{equation}\label{eq:stability}
\| \Gamma_1' \Gamma_1 \| +   \|\Omega\| +  \|\Omega^{-1}\| \lesssim 1,
\end{equation}
hold.

Assume that for some sequence of positive numbers $\{r_n\}$ such that $r_n \to 0$ and $r_n n^{1/2} \to \infty$, the following stochastic equicontinuity and continuity conditions hold:
\begin{eqnarray}\label{eq:se1}
&& \sup_{ \alpha \in \mathcal{A} } \frac{\| \hat \M(\alpha, \hat \eta) - \M(\alpha, \hat \eta) \| + \|  \M(\alpha, \hat \eta) - \M(\alpha,  \eta_0) \|}{r_n + \| \hat \M (\alpha, \hat \eta)\| + \|\M(\alpha,  \eta_0)\| } \to_{\Pr_n} 0,\\
&& \sup_{ \| \alpha - \alpha_0\| \leq r_n} \frac{\| \hat \M(\alpha, \hat \eta) - \M(\alpha, \hat \eta) - \hat \M(\alpha_0, \eta_0) \|}{n^{-1/2} + \| \hat \M (\alpha, \hat \eta)\| + \|\M(\alpha, \eta_0)\| } \to_{\Pr_n} 0.\label{eq:se2}
\end{eqnarray}
Suppose that uniformly for all $\alpha \neq \alpha_0$ such that $\| \alpha - \alpha_0\| \leq r_n \to 0$, the following
conditions on the smoothness of $\M$ and the quality of the estimator $\hat \eta$ hold, as $n \to \infty$:
 \begin{align}\label{eq:derivatives}
 \begin{array}{lll}
 && \| \M(\alpha, \eta_0) - \M(\alpha_0, \eta_0)- \Gamma_1 [\alpha - \alpha_0]\|  \| \alpha - \alpha_0\|^{-1} \to 0,\\
 && \sqrt{n} \| \M(\alpha, \hat \eta) - \M(\alpha, \eta_0)- \partial_{\eta'} \M (\alpha, \eta_0) [\hat \eta - \eta_0]\|  \to_{\Pr_n} 0, \\
 &&   \|  \{ \partial_{\eta'} \M (\alpha, \eta_0) -  \partial_{\eta'} \M (\alpha_0, \eta_0)\} [\hat \eta - \eta_0] \| \|\alpha - \alpha_0\|^{-1}
 \to_{\Pr_n} 0.  \end{array} \end{align}
Finally, as before, we assume that the orthogonality condition
\begin{equation}\label{eq:ortho-e}
\partial_{\eta'} \M (\alpha_0, \eta_0) = 0
\end{equation}
holds.

The above conditions extend the analysis of \citen{pakes:pollard} and \citen{chen:linton:k}, which
in turn extended Huber's (\citeyear{huber}) classical results on Z-estimators.  These conditions allow
for both smooth and non-smooth systems of estimating equations. The identifiability condition  imposed above
is mild and holds for broad classes of identifiable models.  The equicontinuity and smoothness
conditions imposed above require mild smoothness on the function $\M$
and also require that $\hat \eta$ is a good-quality estimator of $\eta_0$. In particular,
these conditions will often require that $\hat \eta$ converges to $\eta_0$ at a faster rate than $n^{-1/4}$
as demonstrated, for example, in the next section.  However, the rate condition alone is not sufficient for adaptivity. We also need
the orthogonality condition (\ref{eq:ortho-e}). In addition, it is required that $\hat \eta \in \mathcal{H}_n$,
where $\mathcal{H}_n$ is a set whose complexity does not grow too quickly with the sample
size, to verify the stochastic equicontinuity condition; see, e.g., \citen{BCFH:Policy} and \citen{BCK-LAD}.
In the next section, we use the sparsity of $\hat \eta$ to control this complexity. Note that 
conditions (\ref{eq:se1})-(\ref{eq:se2})
can be simplified by leaving only $r_n$ and $n^{-1/2}$  in the denominator, though
this simplification would then require imposing compactness on $\mathcal{A}$ even in linear problems.

\begin{proposition}[Valid Inference via Adaptive Estimation after Selection or Regularization]\label{prop:estimation}
Consider a sequence $\{\mathbf{P}_n\}$ of sets
of probability laws such that for each sequence $\{\Pr_n\} \in \{\mathbf{P}_n\}$
conditions (\ref{eq:minimization})-(\ref{eq:ortho-e}) hold.  Then
$$
\sqrt{n}(\hat \alpha -\alpha_0) +  [\Gamma_1'\Gamma_1]^{-1}  \Gamma_1' \sqrt{n} \hat  \M( \alpha_0,  \eta_0) \to_{\Pr_n}0.
$$
In addition, for $V_n := (\Gamma'_1 \Gamma_1)^{-1} \Gamma_1' \Omega \Gamma_1(\Gamma'_1 \Gamma_1)^{-1},$
we have that
$$
\lim_{n \to \infty} \sup_{\Pr \in \mathbf{P}_n} \sup_{R \in \mathcal{R}} |\Pr (  V_n^{-1/2} (\hat \alpha - \alpha_0) \in  R) -  \mathbb{P}( \mathcal{N}(0,I) \in R) |  = 0,
$$
where $\mathcal{R}$ is a collection of all convex sets. Moreover, the result continues to apply if $V_n$ is replaced by a consistent estimator $\hat V_n$ such that
$\hat V_n - V_n \to_{\Pr_n} 0$ under each sequence $\{\Pr_n\}$.   Thus,  $ \mathsf{CR}^l_{1-a} = [
l'\hat \alpha \pm c(1-a/2) (l'\hat V_nl/n)^{1/2}]$ where $c(1-a/2)$ is the $(1-a/2)$-quantile
of $\mathcal{N}(0,1)$ is a uniformly valid confidence set for $l'\alpha_0$:
$$
\lim_{n \to \infty} \sup_{\Pr \in \mathbf{P}_n}  |\Pr (  l'\alpha_0 \in \mathsf{CR}^l_{1-a} ) - (1-a) |  = 0.
$$
\end{proposition}

Note that the above formulation implicitly accommodates weighting options.  Suppose
$\M^o$ and $\hat \M^o$ are the original theoretical and empirical systems of equations,
and let $\Gamma_1^o = \partial_{\alpha'} \M^o (\alpha_0, \eta_0)$ be the original Jacobian.
We could consider $k \times k$ positive-definite weight matrices $\A$ and $\hat{\A}$ such that
\begin{equation}\label{eq:A}
 \|\A^2\| + \| (\A^2)^{-1} \| \lesssim 1,   \quad \|\hat{\A}^2 - \A^2 \| \to_{\Pr_n} 0.
\end{equation}
For example, we may wish to use the optimal weighting matrix $
\A^2 =  \Var(\sqrt{n} \hat{\M}^o(\alpha_0, \eta_0))^{-1}$ which can be estimated by $\hat \A^2$ obtained
using a preliminary estimator $\hat \alpha^o$ resulting from solving the problem with some non-optimal weighting
matrix such as $I$.  We can then simply redefine the system of equations and the Jacobian
according to
\begin{equation}\label{eq:redefine}
\M(\alpha, \eta) =  \A \M^o(\alpha, \eta), \quad  \hat{\M}(\alpha, \eta) = \hat {\A} \hat{\M}^o(\alpha, \eta), \quad \Gamma_1 = \A \Gamma_1^o.
\end{equation}

\bigskip

\begin{proposition}[Adaptive Estimation via Weighted Equations]\label{prop:weighted}
Consider a sequence $\{\mathbf{P}_n\}$ of sets
of probability laws such that for each sequence $\{\Pr_n\} \in \{\mathbf{P}_n\}$
the
conditions of Proposition \ref{prop:estimation} hold for the original pair of
 systems of equations $(\M^o, \hat{\M}^o)$ and that (\ref{eq:A}) holds. Then these
 conditions also hold for  the new pair $(\M, \hat{\M})$ in (\ref{eq:redefine}), so that all the conclusions
of Proposition \ref{prop:estimation} apply to the resulting approximate argmin estimator $\hat \alpha$.  In particular, if
we use $\A^2 =  \Var(\sqrt{n} \hat{\M}^o(\alpha_0, \eta_0))^{-1}$  and $\hat \A^2 - \A^2 \to_{\Pr_n} 0$, then
the large sample variance $V_n$ simplifies to
$
V_n = (\Gamma_1'\Gamma_1)^{-1}.
$

\end{proposition}

\subsection{Inference via Adaptive ``One-Step" Estimation}
We next consider a ``one-step" estimator.  To define the estimator, we start with an initial
estimator $\tilde \alpha$ that satisfies,  for $r_n = o(n^{-1/4})$,
\begin{equation}\label{eq:crude}
\Pr_n(\| \tilde \alpha - \alpha_0 \| \leq r_n) \to 1.
\end{equation}
The one-step estimator $\check\alpha$ then solves a linearized version of  (\ref{eq:minimization}):
\begin{equation}\label{eq:one-step}
\check \alpha = \tilde \alpha -  [\hat \Gamma_1'\hat \Gamma_1]^{-1}  \hat \Gamma_1'  \hat  \M(\tilde \alpha, \hat \eta)
\end{equation}
where $\hat \Gamma_1$ is an estimator of $\Gamma_1$ such that
\begin{equation}\label{eq:consistent gamma}
\Pr_n(\|\hat \Gamma_1 - \Gamma_1\| \leq r_n) \to 1.
\end{equation}

Since the one-step estimator is considerably more crude than the argmin estimator, we need
to impose additional smoothness conditions. Specifically, we suppose that uniformly for all $\alpha \neq \alpha_0$ such that $\| \alpha - \alpha_0\| \leq r_n \to 0$, the following strengthened conditions on stochastic equicontinuity,
smoothness of $\M$ and the quality of the estimator $\hat \eta$ hold, as $n \to \infty$:
 \begin{align}\label{eq:derivatives2}
 \begin{array}{lll}
 &&  n^{1/2}\| \hat \M(\alpha, \hat \eta) - \M(\alpha, \hat \eta) - \hat \M(\alpha_0, \eta_0) \| \to_{\Pr_n} 0,\\
 && \| \M(\alpha, \eta_0) - \M(\alpha_0, \eta_0)- \Gamma_1 [\alpha - \alpha_0]\|  \| \alpha - \alpha_0\|^{-2} \lesssim 1,\\
 && \sqrt{n} \| \M(\alpha, \hat \eta) - \M(\alpha, \eta_0)- \partial_{\eta'} \M (\alpha, \eta_0) [\hat \eta - \eta_0]\|  \to_{\Pr_n} 0, \\
 &&   \sqrt{n} \|  \{ \partial_{\eta'} \M (\alpha, \eta_0) -  \partial_{\eta'} \M (\alpha_0, \eta_0)\} [\hat \eta - \eta_0] \|
 \to_{\Pr_n} 0.  \end{array} \end{align}

\bigskip

\begin{proposition}[Valid Inference via Adaptive One-Step Estimators]\label{prop:one-step}
Consider a sequence $\{\mathbf{P}_n\}$ of sets
of probability laws such that for each sequence $\{\Pr_n\} \in \{\mathbf{P}_n\}$  the
conditions of Proposition \ref{prop:estimation}
 as well as (\ref{eq:crude}), (\ref{eq:consistent gamma}), and (\ref{eq:derivatives2}) hold.  Then
 the one-step estimator $\check \alpha$ defined by (\ref{eq:one-step}) is first order
equivalent to the argmin estimator $\hat \alpha$:
$$
\sqrt{n}(\check \alpha - \hat \alpha) \to_{\Pr_n} 0.
$$
Consequently, all conclusions of Proposition \ref{prop:estimation} apply to
$\check \alpha$ in place of $\hat \alpha$.
\end{proposition}

The one-step estimator requires stronger regularity conditions than the argmin estimator.
Moreover, there is finite-sample evidence (e.g. \cite{BCY-honest}) that in
practical problems the argmin estimator often works much better, since the one-step estimator
typically suffers from higher-order biases.  This problem could be alleviated somewhat
by iterating on the one-step estimator, treating the previous iteration as the ``crude" start $\tilde \alpha$
for the next iteration.

\section{Achieving Orthogonality Using Neyman's Orthogonalization}

Here we describe orthogonalization ideas that go back at least to \citen{Neyman59}; see also \citen{Neyman1979}.
Neyman's idea was to project the score
that identifies the parameter of interest onto the ortho-complement of the tangent
space for the nuisance parameter. This projection underlies semi-parametric
efficiency theory, which is concerned particularly with the case in which $\eta$
is infinite-dimensional, cf. \cite{vdV}.  Here we consider finite-dimensional $\eta$ of high dimension;
for discussion of infinite-dimensional $\eta$ in an approximately sparse setting, see \citen{BCFH:Policy} and \citen{BCK-LAD}.

\subsection{The Classical Likelihood Case} In likelihood settings, the construction of orthogonal equations was proposed by \citen{Neyman59} who used them in construction of his celebrated $C(\alpha)$-statistic.  The $C(\alpha)$-statistic, or the orthogonal score statistic, was first explicitly utilized for testing (and also for setting up estimation) in high-dimensional sparse models in \citen{BCK-LAD} and \citen{BCK-QR},
in the context of quantile regression, and \citen{BCY-honest} in the context of logistic regression and other generalized linear models. More recent uses of $C(\alpha)$-statistics (or close variants) include those by \citen{witten:score}, \citen{hanliu1}, and \citen{hanliu2}.

 Suppose that the (possibly conditional, possibly quasi) log-likelihood function associated with observation
$w_i$ is $\ell(w_i, \alpha, \beta)$, where $\alpha \in \mathcal{A} \subset \mathbb{R}^{d}$
is the target parameter and $\beta \in \mB \subset \mathbb{R}^{p_0}$ is the nuisance parameter.
Under regularity conditions, the true parameter values $\gamma_0=(\alpha_0', \beta_0)'$ obey
\begin{equation}\label{eq: foc lik}
\Ep [\partial_\alpha \ell (w_i, \alpha_0, \beta_0)]=0, \quad \Ep[\partial_{\beta} \ell (w_i, \alpha_0, \beta_0)] =0.
\end{equation}
Now consider the moment function
\begin{equation}
{\M}(\alpha, \eta) = \Ep[\psi(w_i, \alpha, \eta)] ,  \ \ \psi(w_i, \alpha, \eta) =  \partial_\alpha \ell (w_i, \alpha, \beta) -   \mu \partial_{\beta} \ell (w_i, \alpha, \beta).
\end{equation}
Here the nuisance parameter is $$  \eta= (\beta', \textrm{vec}(\mu)')'  \in \mB \times \mathcal{D} \subset \mathbb{R}^{p}, \quad p=p_0 + dp_0,$$
where $\mu$ is the $d \times p_0$ \textit{orthogonalization} parameter matrix whose
true value  $\mu_0$ solves the equation:
\begin{equation}\label{eq: ortho gono}
J_{\alpha\beta}  - \mu J_{\beta \beta} =0  \ (\text{ i.e., }  \mu_0 =  J_{\alpha\beta}J^{-1}_{\beta \beta}),
\end{equation}
where, for $\gamma := (\alpha', \beta')'$ and $\gamma_0 := (\alpha_0', \beta_0')'$,
\begin{eqnarray*}
J := - \partial_{\gamma'} \Ep [   \partial_{\gamma}  \ell(w_i, \gamma)\ ] \vert_{\gamma = \gamma_0}                                                                                  &=: &  \left(
                                                                                   \begin{array}{cc}
                                                                                     J_{\alpha \alpha} & J_{\alpha \beta} \\
                                                                                     J_{\beta \alpha} &  J_{\beta \beta} \\
                                                                                   \end{array}
                                                                                 \right).
\end{eqnarray*}
Note that $\mu_0$ not only creates the necessary orthogonality but also creates
\begin{itemize}
\item the \textit{optimal score} (in statistical language)
\item or, equivalently, the \textit{optimal instrument/moment} (in econometric language)\footnote{The connection
between optimal instruments/moments and likelihood/score has been elucidated by the fundamental work of \cite{chamberlain}.}
\end{itemize}
for inference about $\alpha_0$. 


Provided $\mu_0$ is well-defined, we have by (\ref{eq: foc lik}) that
$${\M}(\alpha_0, \eta_0) = 0.$$ Moreover, the function $\M$ has the desired orthogonality property:\begin{equation}\label{eq: ortho lik}
\partial_{\eta'} {\M}(\alpha_0, \eta_0) = \Big [J_{\alpha \beta} - \mu_0 J_{\beta \beta}; \   F \Ep[\partial_\beta  \ell (w_i, \alpha_0, \beta_0)]  \Big ] =0,
\end{equation}
where $F$ is a tensor operator, such that $F x = \partial \mu x/ \partial \mathrm{vec(\mu)}' \mid_{\mu = \mu_0}$
is a $d \times (d p_0)$ matrix for any vector $x$ in $\mathbb{R}^{p_0}$.  Note that the orthogonality property holds for Neyman's construction even if the likelihood is misspecified. That is, $\ell(w_i, \gamma_0)$
may be a quasi-likelihood, and the data need not be i.i.d. and may, for example,
exhibit complex dependence over $i$.

An alternative way to define $\mu_0$ arises by considering
that,  under the correct specification and sufficient regularity, the information matrix equality holds  and yields
\begin{eqnarray*}
J = J^0 & := & \Ep [ \partial_\gamma \ell(w_i, \gamma) \partial_\gamma \ell(w_i, \gamma) '] \vert_{\gamma = \gamma_0}  \\
& = & \left .  \left(  \begin{array}{cc}
                                                                                     \Ep[ \partial_\alpha \ell(w_i, \gamma) \partial_{\alpha} \ell(w_i, \gamma)']  &  \Ep[ \partial_\alpha \ell(w_i, \gamma) \partial_\beta \ell(w_i, \gamma)'] \\
                                                                                      \Ep[ \partial_\beta \ell(w_i, \gamma) \partial_\alpha \ell(w_i, \gamma)']&     \Ep[ \partial_\beta \ell(w_i, \gamma) \partial_\beta \ell(w_i, \gamma)']\\
                                                                                   \end{array}
                                                                                 \right) \right |_{\gamma = \gamma_0}, \\
                                                                                 &=: &  \left(
                                                                                   \begin{array}{cc}
                                                                                     J^0_{\alpha \alpha} & J^0_{ \alpha \beta} \\
                                                                                     J^{0}_{\beta \alpha} &  J^0_{\beta \beta} \\
                                                                                   \end{array}
                                                                                 \right).
\end{eqnarray*}
Hence define  $\mu^*_0 = J^0_{\alpha\beta} J^{0-1}_{\beta \beta}$ as the population \textit{projection coefficient} of the score for the main parameter $ \partial_\alpha \ell(w_i, \gamma_0)$  on the score for the nuisance parameter $\partial_\beta \ell(w_i, \gamma_0)$:
\begin{equation}\label{projection}
\partial_\alpha \ell(w_i, \gamma_0) = \mu^*_0 \partial_\beta \ell(w_i, \gamma_0) +  \varrho, \ \ \Ep[ \varrho \partial_\beta \ell(w_i, \gamma_0)']=0.
\end{equation}
 We can see this construction as  the non-linear version of  Frisch-Waugh's ``partialling out" from the linear regression model.
It is important to note that under misspecification the information matrix equality generally does not hold, and this projection approach does not provide valid orthogonalization.

\begin{lemma}[Neyman's orthogonalization for (quasi-) likelihood scores]
Suppose that for each $\gamma = (\alpha, \beta) \in \mathcal{A}\times\mathcal{B}$,  the derivative $\partial_\gamma \ell(w_i, \gamma)$ exists and is continuous at  $\gamma$ with probability one, and obeys the dominance condition $\Ep \sup_{\gamma \in \mathcal{A}\times\mathcal{B}}\| \partial_\gamma \ell(w_i, \gamma) \|^2< \infty$. Suppose that
condition (\ref{eq: foc lik}) holds for some (quasi-) true value $(\alpha_0, \beta_0)$. Then, (i) if $J$  exists and is finite and $J_{\beta \beta}$ is invertible, then the orthogonality condition (\ref{eq: ortho lik}) holds; (ii)  if 
the information matrix equality holds, namely $J= J^0$, then the orthogonality condition (\ref{eq: ortho lik}) holds for the projection parameter $\mu^*_0$ in place of the orthogonalization parameter matrix $\mu_0$.
\end{lemma}

 The  claim follows immediately from the computations above.
\bigskip

With the formulations given above Neyman's $C(\alpha)$-statistic takes the form
$$
C(\alpha) = \| S(\alpha)\|_2^2, \quad S (\alpha) = \hat \Omega^{-1/2} (\alpha, \hat \eta) \sqrt{n} \hat{\M}(\alpha, \hat \eta),
$$
where $\hat{\M}(\alpha, \hat \eta) = \En [\psi(w_i, \alpha, \hat \eta)]$ as before,
$ \Omega (\alpha,  \eta_0) = \mathrm{Var} ( \sqrt{n} \hat{\M}(\alpha, \eta_0))$,
and $\hat \Omega (\alpha, \hat \eta)$ and $\hat \eta$ are suitable estimators
based on sparsity or other structured assumptions. The estimator is then
$$
\hat \alpha =  \arg\inf_{\alpha \in \AA} C(\alpha) = \arg\inf_{\alpha \in \AA} \| \sqrt{n}\hat{\M}(\alpha, \hat \eta) \|,
$$
provided that $\hat \Omega (\alpha, \hat \eta)$ is positive definite for each $\alpha \in \AA$.  If the conditions
of Section 2 hold, we have that
\begin{equation}\label{mle:result}
 C(\alpha)  \leadsto \chi^2(d), \quad V_n^{-1/2}\sqrt{n}(\hat \alpha - \alpha_0) \leadsto \mathcal{N}(0,I),
 \end{equation}
where $V_n = \Gamma_1^{-1} \Omega (\alpha_0,  \eta_0)\Gamma_1^{-1}$ and $\Gamma_1 =
J_{\alpha \alpha} - \mu_0 J'_{\alpha \beta}$. Under the correct specification and i.i.d. sampling,
the variance matrix $V_n$ further reduces to the optimal
variance $$\Gamma_1^{-1}= (J_{\alpha \alpha} - J_{\alpha \beta} J^{-1}_{\beta \beta }J'_{\alpha \beta})^{-1},$$ of the
first $d$ components of the maximum likelihood estimator in a Gaussian shift experiment with observation $Z \sim \mathcal{N}(h, J^{-1}_0)$. Likewise, the result (\ref{mle:result}) also holds for
the one-step estimator $\check \alpha$  of Section 2 in place of $\hat \alpha$ as long as the  conditions in Section 2 hold.

Provided that sparsity or its generalizations are plausible assumptions to make regarding
$\eta_0$, the formulations above naturally lend themselves to sparse estimation.    For example, \citen{BCY-honest} used
penalized and post-penalized maximum likelihood to estimate $\beta_0$, and used the information matrix
equality to estimate the orthogonalization parameter matrix $\mu^*_0$ by using Lasso or Post-Lasso estimation
of the projection equation (\ref{projection}).    It is also possible to estimate $\mu_0$ directly by
finding approximate sparse solutions to the empirical analog of the system of equations
$
J_{\alpha \beta} - \mu J_{\beta \beta}= 0$ using $\ell_1$-penalized estimation, as, e.g., in \cite{vdGBRD:AsymptoticConfidenceSets}, or post-$\ell_1$-penalized estimation.

\subsection{Achieving Orthogonality in GMM Problems }
 Here we consider $\gamma_0 = (\alpha_0', \beta_0')'$ that solve
 the system of equations: $$
 \Ep [m( w_i, \alpha_0, \beta_0)] =  0,
 $$
 where $m: \mathcal{W} \times \mathcal{A} \times \mB  \mapsto \mathbb{R}^k$,
 $\mathcal{A}\times\mB$ is a convex subset of $\mathbb{R}^{d} \times \mathbb{R}^{p_0}$,
 and $k \geq d+ p_0$ is the number of moments.  The orthogonal moment equation is
 \begin{equation}\label{eq: GMM1}
{\M}(\alpha, \eta) = \Ep[\psi(w_i, \alpha, \eta)] ,  \ \ \psi(w_i, \alpha, \eta) =
\mu m(w_i, \alpha, \beta).
 \end{equation}
The nuisance parameter is $$  \eta= (\beta', \textrm{vec}(\mu)')'  \in \mB \times \mathcal{D} \subset \mathbb{R}^{p}, \quad p=p_0 + dk,$$
where $\mu$ is the $d \times k$ orthogonalization parameter matrix. The ``true value" of $\mu$ is $$
\mu_0=(G_{\alpha}' \Omega_m^{-1} -  G_{\alpha}' \Omega_m^{-1} G_\beta (G_\beta' \Omega_m^{-1} G_\beta)^{-1}  G_\beta '\Omega^{-1}_m  ),
$$
where, for $\gamma = (\alpha', \beta')'$ and $\gamma_0 = (\alpha_0', \beta_0')'$,
$$
G_{\gamma} =  \partial_{\gamma'} \Ep[m(w_i, \alpha, \beta)] \Big |_{\gamma= \gamma_0}= \Big [\partial_{\alpha'} \Ep[m(w_i, \alpha, \beta)],    \partial_{\beta'} \Ep[m(w_i, \alpha, \beta)] \Big ] \Big |_{\gamma= \gamma_0}= : \Big [G_{\alpha}, G_{\beta} \Big],
$$
and
$$
\Omega_m = \Var( \sqrt{n} \En[m(w_i, \alpha_0, \beta_0) ] ).
$$
As before, we can interpret $\mu_0$ as an operator creating orthogonality while building
\begin{itemize}
\item the \textit{optimal instrument/moment} (in econometric language),
\item  or, equivalently,
the \textit{optimal score} function (in statistical language).\footnote{Cf.~previous footnote.}

\end{itemize}
The resulting moment function has the required orthogonality property; namely, the first derivative with respect to the nuisance
parameter when evaluated at the true parameter values is zero:
\begin{equation}\label{GMM:ortho}
\partial_{\eta'}{\M}(\alpha_0, \eta)\vert_{\eta = \eta_0} =  [\mu_0 G_\beta,  F \Ep[m(w_i, \alpha_0, \beta_0)]] = 0,
\end{equation}
 where $F$ is a tensor operator, such that $F x = \partial \mu x/ \partial \mathrm{vec(\mu)}' \mid_{\mu = \mu_0}$
is a $d \times (dk)$ matrix for any vector $x$ in $\mathbb{R}^{k}$. 

Estimation and inference on $\alpha_0$ can be based on the empirical analog of  (\ref{eq: GMM1}):
$$
\hat{\M}(\alpha, \hat \eta) =  \En[ \psi (w_i, \alpha, \hat \eta) ],
$$
where $\hat \eta$ is a post-selection or other regularized estimator of $\eta_0$.
Note that the previous framework of (quasi)-likelihood is incorporated as a special case
with $$m(w_i, \alpha, \beta) = [ \partial_{\alpha} \ell(w_i, \alpha)', \partial_{\beta} \ell(w_i, \beta)']'.$$

With the formulations above, Neyman's $C(\alpha)$-statistic takes the form:
$$
C(\alpha) = \| S(\alpha)\|_2^2, \quad S (\alpha) = \hat \Omega^{-1/2} (\alpha, \hat \eta) \sqrt{n} \hat{\M}(\alpha, \hat \eta),
$$
where $\hat{\M}(\alpha, \hat \eta) = \En [\psi(w_i, \alpha, \hat \eta)]$ as before,
$ \Omega (\alpha,  \eta_0) = \mathrm{Var} ( \sqrt{n} \hat{\M}(\alpha, \eta_0))$,
and $\hat \Omega (\alpha, \hat \eta)$ and $\hat \eta$ are suitable estimators
based on structured assumptions.  The estimator is then
$$
\hat \alpha =  \arg\inf_{\alpha \in \AA} C(\alpha) = \arg\inf_{\alpha \in \AA} \| \sqrt{n}\hat{\M}(\alpha, \hat \eta) \|,
$$
provided that $\hat \Omega (\alpha, \hat \eta)$ is positive definite for each $\alpha \in \AA$.  If the high-level conditions
of Section 2 hold, we have that
 \begin{equation}\label{GMM:result}
 C(\alpha)  \leadsto_{\Pr_n} \chi^2(d), \quad V_n^{-1/2}\sqrt{n}(\hat \alpha - \alpha) \leadsto_{\Pr_n} \mathcal{N}(0,I),
 \end{equation}
where $V_n =(\Gamma_1')^{-1} \Omega (\alpha_0,  \eta_0)(\Gamma_1)^{-1}$ coincides with the optimal variance
for GMM; here $\Gamma_1 = \mu_0 G_{\alpha}$. Likewise, the same result (\ref{GMM:result}) holds for
the one-step estimator $\check \alpha$  of Section 2 in place of $\hat \alpha$ as long as the conditions in Section 2 hold.
 In particular, the variance $V_n$
corresponds to the variance of the first $d$ components
of the maximum likelihood estimator in the normal shift experiment with the observation
$Z \sim \mathcal{N}(h, (G_\gamma'\Omega_m^{-1}G_\gamma)^{-1})$.

The above
is a generic outline of the properties that are expected for inference using orthogonalized GMM equations under structured assumptions. The problem of inference in GMM under sparsity is a very delicate matter due to the complex form of the orthogonalization parameters.
One approach to the problem is developed in \cite{GMM:sparse}.

\section{Achieving Adaptivity In Affine-Quadratic Models via Approximate Sparsity}
Here we take orthogonality as given
and explain how we can use approximate sparsity to achieve the adaptivity property (\ref{eq:adaptivity}).

\subsection{The Affine-Quadratic Model} We analyze the case in which $\hat \M$ and $\M$ are \textit{affine} in $\alpha$ and \textit{affine-quadratic} in $\eta$.  Specifically, we suppose that for
all $\alpha$
$$
\hat \M (\alpha, \eta) = \hat \Gamma_1(\eta)\alpha +  \hat \Gamma_2(\eta), \quad
 \M (\alpha, \eta) = \Gamma_1(\eta)\alpha + \Gamma_2(\eta),
$$
where the orthogonality condition holds,
\begin{align*}
\partial_{\eta'} \M (\alpha_0, \eta_0) = 0,
\end{align*}
and $\eta \mapsto \hat \Gamma_j (\eta)$ and  $\eta \mapsto \Gamma_j (\eta)$ are affine-quadratic in $\eta$ for $j=1$ and $j=2$. That is, we will have that all second-order derivatives of $\hat \Gamma_j (\eta)$ and $\Gamma_j (\eta)$ for $j = 1$ and $j=2$ are constant over the convex parameter
space $\mathcal{H}$ for $\eta$.

This setting is both useful, including most widely used linear models as a special case, and pedagogical, permitting simple illustration of the key issues that arise in treating the general problem.  The derivations
given below easily generalize to more complicated models, but we defer the details to the interested reader.

The estimator in this case is
\begin{equation}\label{eq:linear estimator}
\hat \alpha = \arg\min_{\alpha \in \mathbb{R}^d} \|\hat \M (\alpha, \hat \eta) \|^2= -  [\hat \Gamma_1(\hat \eta)'\hat \Gamma_1(\hat \eta)]^{-1} \hat \Gamma_1(\hat \eta)' \hat  \Gamma_2(\hat \eta),
\end{equation}
provided the inverse is well-defined.  It follows that
\begin{equation}\label{eq:linear estimator represent}
\sqrt{n} (\hat \alpha - \alpha_0) =  - [\hat \Gamma_1(\hat \eta)'\hat \Gamma_1(\hat \eta)]^{-1}  \hat \Gamma_1(\hat \eta)'  \sqrt{n} \hat \M(\alpha_0, \hat \eta).
\end{equation}
This estimator is adaptive if, for $\Gamma_1:= \Gamma_1(\eta_0)$,
$$
\sqrt{n} (\hat \alpha - \alpha_0) + [\Gamma_1' \Gamma_1]^{-1} \Gamma_1'    \sqrt{n}\hat \M(\alpha_0, \eta_0) \to_{\Pr_n} 0,
$$
which occurs under  the conditions in (\ref{eq:normal score}) and  (\ref{eq:stability})  if
\begin{equation}\label{eq:adaptivity estimation linear}
\sqrt{n}(\hat \M(\alpha_0, \hat \eta) - \hat \M(\alpha_0, \eta_0)) \to_{\Pr_n} 0,  \quad  \hat \Gamma_1(\hat \eta) - \Gamma_1(\eta_0)
\to_{\Pr_n} 0.
\end{equation}
Therefore, the problem of the adaptivity of the estimator is directly connected to the problem of the adaptivity of testing hypotheses
about $\alpha_0$.

\bigskip

\begin{lemma}[Adaptive Testing and Estimation in Affine-Quadratic Models] \label{prop:affine}
Consider a sequence $\{\mathbf{P}_n\}$ of sets  of probability laws such that for each sequence $\{\Pr_n\} \in \{\mathbf{P}_n\}$, conditions stated in the first paragraph of Section 4.1, condition (\ref{eq:adaptivity estimation linear}), the asymptotic normality condition
 (\ref{eq:normal score}), the stability condition (\ref{eq:stability}), and condition (\ref{eq:variance}) hold. Then
 all the conditions of Propositions 1 and 2 hold. Moreover, the conclusions of Proposition 1 hold, and the conclusions of Proposition \ref{prop:estimation} hold for the estimator $\hat \alpha$ in   (\ref{eq:linear estimator}).
\end{lemma}

\subsection{Adaptivity for Testing via Approximate Sparsity}
Assuming the
orthogonality condition holds, we follow \cite{BCCH12}  in  using approximate
sparsity to achieve the adaptivity property (\ref{eq:adaptivity}) for the testing problem in the
affine-quadratic models.

We can expand each element $\hat{\M}_j$ of $\hat{\M}= (\hat \M_j )_{j=1}^k$
as follows:
\begin{equation}\label{eq:decompose}
\sqrt{n}(\hat{\M}_j(\alpha_0, \hat \eta) - \hat{\M}_j(\alpha_0, \eta_0)) =  T_{1,j} + T_{2,j} +T_{3,j},
\end{equation}
where
\begin{equation}\label{eq:terms}
\begin{array}{lll}
& &  T_{1,j}:=  \sqrt{n}\partial_{\eta} {\M}_j(\alpha_0, \eta_0) '(\hat \eta - \eta_0), \\
&& T_{2,j}:= \sqrt{n} (\partial_\eta \hat{\M}_j(\alpha_0, \eta_0) -  \partial_\eta  {\M}_j(\alpha_0, \eta_0) ) '(\hat \eta - \eta_0), \\
& &  T_{3,j}:=  \sqrt{n} 2^{-1}   (\hat \eta - \eta_0)' \partial_\eta  \partial_{\eta'}\hat{\M}_j(\alpha_0) (\hat \eta - \eta_0).
\end{array}
\end{equation}
The term $T_{1,j}$ vanishes precisely because of orthogonality, i.e. $$T_{1,j}=0.$$ However, terms $T_{2,j}$ and $T_{3,j}$
need not vanish.  In order to show that they are asymptotically negligible, we need to impose further structure on the problem.

\textbf{Structure 1: Exact Sparsity.} We first consider the case of using an exact sparsity structure where $\|\eta_0\|_0 \leq s$ and $s=s_n \geq 1$ can depend on $n$.  We then use
estimators $\hat \eta$ that exploit the sparsity structure.

Suppose that the following bounds hold with probability $1 - o(1)$ under $\Pr_n$:
\begin{equation}\label{eq:sparse}
\begin{array}{c}
 \|\hat \eta\|_0 \lesssim s, \quad   \|\eta_0 \|_0 \leq s,  \\
 \| \hat \eta- \eta_0\|_2 \lesssim \sqrt{(s/n) \log (pn)}, \quad   \|\hat \eta - \eta_0\|_1
\lesssim  \sqrt{(s^2/n) \log (pn)}.
\end{array}
\end{equation}
These conditions are typical performance bounds which hold for many sparsity-based
estimators such as Lasso, post-Lasso, and their extensions.

We suppose further that the moderate deviation bound
\begin{equation}\label{eq:SNMD}
\bar T_{2,j}=\| \sqrt{n}  (\partial_{\eta'} \hat{\M}_j(\alpha_0, \eta_0) -  \partial_{\eta'} {\M}_j(\alpha_0, \eta_0) )\|_{\infty} \lesssim_{\Pr_n} \sqrt{\log (p n)}
\end{equation}
holds and that the sparse norm of the second-derivative matrix is bounded:
\begin{equation}\label{eq:BDD}
\bar T_{3,j}= \|  \partial_\eta \partial_{\eta'}  \hat{\M}_j ( \alpha_0) \|_{\mathsf{sp}(\ell_n s)} \lesssim_{\Pr_n} 1
\end{equation}
where  $\ell_n \to \infty$ but $\ell_n =o( \log n)$.

Following \cite{BCCH12}, we can verify condition (\ref{eq:SNMD}) using the moderate deviation theory for self-normalized sums (e.g., \cite{jing:etal}), which allows us to avoid making highly restrictive subgaussian or gaussian tail assumptions. Likewise, following \cite{BCCH12}, we can verify the second condition using laws of large numbers for large
matrices acting on sparse vectors as in \cite{rudelson:vershynin} and \citen{RudelsonZhou2011}; see Lemma 7.  Indeed, condition (\ref{eq:BDD}) holds if
$$
\begin{array}{ll}
\|  \partial_\eta \partial_{\eta'}  \hat{\M}_j ( \alpha_0) -\partial_\eta \partial_{\eta'}   {\M}_j(\alpha_0)   \|_{\mathsf{sp}(\ell_n s)} \to_{\Pr_n} 0,  \
\|\partial_\eta \partial_{\eta'}   {\M}_j ( \alpha_0)   \|_{\mathsf{sp}(\ell_n s)} \lesssim 1.
\end{array}
$$

The above analysis immediately implies the following elementary result.

\bigskip

\begin{lemma}[Elementary Adaptivity for Testing via Sparsity]\label{lemma:ad testing via sparsity} Let $\{ \Pr_n\}$ be a sequence of probability laws. Assume (i)  $\eta \mapsto \hat \M(\alpha_0, \eta)$ and $\eta \mapsto \M(\alpha_0, \eta)$
 are affine-quadratic in $\eta$ and the orthogonality condition holds, (ii) that the conditions on sparsity and the quality of estimation (\ref{eq:sparse}) hold, and the sparsity index obeys
\begin{equation}\label{strong}
s^2 \log (pn)^2/n \to 0,
\end{equation}
(iii) that the moderate deviation bound (\ref{eq:SNMD}) holds, and (iv) the sparse norm of the second
derivatives matrix is bounded as in (\ref{eq:BDD}). Then the adaptivity condition (\ref{eq:adaptivity}) holds for the sequence  $\{ \Pr_n\}$.
\end{lemma}

We note that  (\ref{strong}) requires that the true value of the nuisance parameter is sufficiently sparse,
which we can relax in some special cases to the requirement $s \log (pn)^{c}/n \to 0$, for some constant $c$, by using
sample-splitting techniques; see \cite{BCCH12}. However, this requirement seems unavoidable in general.

Proof. We note above that $T_{1,j} = 0$ by orthogonality. Under (\ref{eq:sparse})-(\ref{eq:SNMD})  if   $s^2 \log (pn)^2/n \to 0$, then $T_{2,j}$ vanishes in probability, as by H\"{o}lder's inequality,
$$
T_{2,j} \leq \bar T_{2,j}    \|\hat \eta - \eta_0\|_1  \lesssim_{\Pr_n}  \sqrt{s^2 \log (pn)^2/n} \to_{\Pr_n} 0.
$$
Also, if $s^2 \log (pn)^2/n \to 0$, then $T_{3,j}$ vanishes in probability, since by H\"{o}lder's inequality
and for sufficiently large $n$,
$$
T_{3,j} \leq \bar T_{3,j} \|\hat \eta - \eta_0\|^2 \lesssim_{\Pr_n} \sqrt{n} s \log (pn)/n \to_{\Pr_n} 0.
$$
The conclusion follows from (\ref{eq:decompose}). \qed

\textbf{Structure 2.  Approximate Sparsity.} Following \cite{BCCH12},
we next consider an approximate sparsity structure.  Approximate sparsity imposes that,
given a constant $c>0$,  we can decompose $\eta_0$ into a
sparse component $\eta^m_r$ and a ``small" non-sparse component $\eta^r$:
\begin{align}\label{eq:app sparsity}
\begin{array}{c}
\eta_0 = \eta^m_{0} +  \eta^r_{0},  \ \spt(\eta^m_{0}) \cap \spt(\eta^r_{0}) = \emptyset,  \\
  \| \eta_{0}^m\|_0 \leq s,  \  \|\eta_{0}^r\|_2 \leq c \sqrt{s/n},    \ \|\eta^r_{0}\|_1 \leq c\sqrt{s^2/n}.
  \end{array}
\end{align}
This condition allows for much more realistic and
richer models than can be accommodated under exact sparsity.  For example,  
$\eta_0$ needs not have \textit{any} zero components at all under approximate sparsity.
In Section 5, we provide an example in which (\ref{eq:app sparsity}) arises from a
more primitive condition that the absolute values $\{|\eta_{0j}|, j =1,...,p\}$, sorted in decreasing order,
decay at a polynomial speed with respect to $j$.

 Suppose that we have an estimator $\hat \eta$ such that
 with probability $1 - o(1)$ under $\Pr_n$ the following bounds hold:
\begin{equation}\label{eq:sparse2}
\begin{array}{c}
 \|\hat \eta\|_0 \lesssim s,
 \| \hat \eta- \eta^m_0\|_2 \lesssim \sqrt{(s/n) \log (pn)}, \quad   \|\hat \eta - \eta^m_0\|_1
\lesssim  \sqrt{(s^2/n) \log (pn)}.
\end{array}
\end{equation}
This condition is again a standard performance bound expected to hold for sparsity-based
estimators under approximate sparsity conditions; see \cite{BCCH12}. Note that by the approximate sparsity condition, we also have that, with probability $1- o(1)$ under $\Pr_n$,
\begin{equation}\label{eq:sparse overall performance}
 \| \hat \eta- \eta_0\|_2 \lesssim \sqrt{(s/n) \log (pn)}, \quad   \|\hat \eta - \eta_0\|_1
\lesssim  \sqrt{(s^2/n) \log (pn)}.
\end{equation}

We can employ the same moderate deviation and bounded sparse norm
conditions as in the previous subsection. In addition, we require the pointwise norm of the second-derivatives matrix to be bounded.  Specifically,
for any deterministic vector $a \neq 0$, we require
\begin{equation}\label{eq:pointwise BDD}
\ \|\partial_\eta  \partial_{\eta'}\hat{\M}_j(\alpha_0)\|_{{\sf pw}(a)} \lesssim_{\Pr_n} 1.
\end{equation}
This condition can be easily verified using ordinary laws of large numbers.

\bigskip

\begin{lemma}[Elementary Adaptivity for Testing via Approximate Sparsity]\label{lemma:adapt approx sparsity} Let $\{ \Pr_n\}$ be a sequence of probability laws. Assume (i)  $\eta \mapsto \hat \M(\alpha_0, \eta)$ and $\eta \mapsto \M(\alpha_0, \eta)$
 are affine-quadratic in $\eta$ and the orthogonality condition holds,  (ii) that the conditions on approximate sparsity (\ref{eq:app sparsity}) and the quality of estimation (\ref{eq:sparse2}) hold, and the sparsity index obeys
$$
s^2 \log (pn)^2/n \to 0,
$$
(iii) that the moderate deviation bound (\ref{eq:SNMD}) holds,  (iv) the sparse norm of the second
derivatives matrix is bounded as in (\ref{eq:BDD}), and (v) the pointwise norm of the second derivative
matrix is bounded as in (\ref{eq:pointwise BDD}). Then the adaptivity condition (\ref{eq:adaptivity}) holds:
$$\sqrt{n}(\hat{\M}(\alpha_0, \hat \eta) - \hat{\M}(\alpha_0, \eta_0))    \to_{\Pr_n} 0.$$
\end{lemma}

\subsection{Adaptivity for Estimation via Approximate Sparsity} We work with the approximate sparsity setup and the affine-quadratic model introduced in the previous subsections.

In addition to the previous assumptions, we  impose the following conditions on the components
$\partial_\eta  {\G}_{1,ml}$ of $\partial_\eta  {\G}_{1}$, where $m=1,...,k$ and $l=1,...,d$,. First, we need the following
deviation and boundedness condition: For each $m$ and $l$,
\begin{equation}\label{eq:SNMD estimation}
 \|\partial_\eta \hat{\G}_{1,ml}(\eta_0) -  \partial_\eta  {\G}_{1,ml}(\eta_0)\|_{\infty}
 \lesssim_{\Pr_n} 1,  \quad \|\partial_{\eta} {\G}_{1, ml}(\eta_0) \|_{\infty} \lesssim 1.
\end{equation}
Second, we require the sparse and pointwise norms of the following second-derivative
matrices be stochastically bounded: For each $m$ and $l$,
\begin{equation}\label{eq:Gamma BDD}
\|\partial_\eta  \partial_{\eta'}\hat{\G}_{1, ml}\|_{\mathsf{sp}(\ell_n s)}  + \|\partial_\eta  \partial_{\eta'}\hat{\G}_{1, ml}\|_{\mathsf{pw}(a)} \lesssim_{\Pr_n} 1,
\end{equation}
where $a \neq 0$ is any deterministic vector. Both of these conditions are mild. They can be verified using self-normalized
moderate deviation theorems and by using laws of large numbers for matrices
as discussed in the previous subsection.

\begin{lemma}[Elementary Adaptivity for Estimation via Approximate Sparsity]\label{lemma:est adapt approx sparsity}
Consider a sequence $\{\Pr_n\}$ for which the conditions of the previous lemma hold. In addition assume that the deviation
 bound (\ref{eq:SNMD estimation}) holds and  the sparse norm and pointwise norms of the second
derivatives matrices are stochastically bounded as in (\ref{eq:Gamma BDD}).
Then the adaptivity condition (\ref{eq:adaptivity estimation linear}) holds for the testing and estimation problem
in the affine-quadratic model.\end{lemma}

\section{Analysis of the IV Model with Very Many Control and Instrumental Variables}
Note that in the following
we write $w \perp v$ to denote $\Cov(w,v) = 0$.

Consider the linear instrumental variable model with response variable:
\begin{align}\label{ystructure}
\begin{array}{lll}
 y_i =   d_i'\alpha_0 + x_i' \beta_0 + \varepsilon_i, &  \Ep[\varepsilon_i]= 0, &  \varepsilon_i \perp (z_i, x_i),
 \end{array}
\end{align}
where $y_i$ is the response variable, $d_i = (d_{ik})_{k=1}^{p^d}$ is a $p^d$-vector of endogenous variables,
such that
\begin{align}
\label{dstructure}
\begin{array}{lll}
d_{i1} \ =  x_{i}'\gamma_{01} \ + z_{i}'\delta_{01} \ + u_{i1}, &  \ \Ep[u_{i1}]= 0,  & \ u_{i1}  \perp (z_i, x_i), \\
\vdots & \vdots & \vdots \\
d_{ip^d} =  x_{i}'\gamma_{0p^d} + z_{i}'\delta_{0p^d} + u_{ip^d}, &   \Ep[u_{ip^d}]= 0, & u_{ip^d}  \perp (z_i, x_i). \\
\end{array} \end{align}
Here $x_i=(x_{ij})_{j=1}^{p^x}$ is a $p^x$-vector of exogenous control variables, including a constant, and $z_i
=(z_i)_{i=1}^{p^z}$ is a $p^z$-vector of instrumental variables.  
We will have $n$ \textit{i.i.d.} draws of $w_i = (y_i, d'_i, x'_i, z'_i)'$
 obeying this system of equations.  We also assume that $\Var(w_i)$ is finite throughout so that the model is well defined.

  The parameter value $\alpha_0$ is our target.  We allow $p^x= p_n^x \gg n$ and $p^z = p_n^z \gg n$,
but we maintain that $p^d$ is fixed in our analysis.    This model includes the case of many instruments and small number of controls considered by \cite{BCCH12} as a special case, and the analysis readily accommodates the case of many controls and no instruments -- i.e. the linear
regression model -- considered by \cite{BCH2011:InferenceGauss,BelloniChernozhukovHansen2011} and \cite{ZhangZhang:CI}. For the latter, we simply set $p_n^z = 0$ and impose the additional condition $\varepsilon_i \perp u_i$
for $u_i = (u_{ij})_{j=1}^{p_d}$, which together with $\varepsilon_i \perp x_i$ implies that $\varepsilon_i \perp d_i$. We also note that the condition $\varepsilon_i \perp x_i, z_i$ is weaker than the condition $\Ep[\varepsilon_i|x_i,z_i]=0$, which allows for some misspecification of the model.

We may have that $z_i$ and $x_i$ are correlated so that $z_i$ are valid instruments only after controlling for $x_i$; specifically, we let $z_i = \Pi x_i + \zeta_i,$ for $\Pi$ a $p_n^z \times p_n^x$ matrix and $\zeta_i$ a $p_n^z$-vector of unobservables with $x_i \perp \zeta_{i}$.  Substituting this expression for $z_i$ as a function of $x_i$ into  (\ref{ystructure}) gives a system for $y_i$ and $d_i$ that depends only on $x_i$:
\begin{align}
\begin{array}{lll}
 y_i \ = x_{i}'\theta_0 \ + \ \rho^y_i, & \Ep[\rho^y_i]\ = 0, &  \rho^y_i \perp x_i,  \\
\ & \ &  \\
d_{i1} \ \ = x_{i}'\vartheta_{01} + \rho^d_{i1}, &  \Ep[\rho^d_{i1}]= 0, &  \rho^d_{i1} \perp x_i, \\
\vdots & \vdots & \vdots \\
d_{ip^{d}} = x_{i}'\vartheta_{0p^d} + \rho^d_{ip^d}, &  \Ep[\rho^d_{ip^d}]= 0, &  \rho^d_{ip^d} \perp x_i. \\
\end{array}
\end{align}


Because the dimension $p=p_n$ of $$\eta_0 = ({\theta}'_0, ({\vartheta}_{0k}', \gamma_{0k}', \delta_{0k}')_{k=1}^{p^d} )'$$ may be larger than $n$, informative estimation and inference about $\alpha_0$ is impossible without imposing restrictions on $\eta_0$.

To state our assumptions, we fix a collection of positive constants $(\sa, \sA, \sc, \sC)$, where $\sa>1$, and a sequence of constants $\delta_n \searrow 0$ and $\ell_n \nearrow \infty$.   These constants will not vary with $\Pr$, but rather we will work with collections of $\Pr$ defined  by these constants.

\textsc{Condition AS.1}  \textit{ We assume that
$\eta_0$ is approximately sparse, namely that the decreasing rearrangement $(|\eta_0|^*_{j})_{j=1}^p $
of absolute values of coefficients  $(|\eta_{0j}|)_{j=1}^p$ obeys}
\begin{align}
| \eta_0|^*_j \leq \mathsf{A} j^{-\mathsf{a}}, \ \  \mathsf{a} > 1,  \ \ j = 1,...,p.
\end{align}

Given this assumption  we can decompose $\eta_0$ into a
sparse component $\eta^m_0$ and small non-sparse component $ \eta^r_{0}$:
\begin{align}\label{decomposition}
\begin{array}{c}
\eta_0 = \eta^m_{0} +  \eta^r_{0},  \ \spt(\eta^m_{0}) \cap \spt(\eta^r_{0}) = \emptyset,  \\
  \| \eta^m_{0}\|_0 \leq s,  \  \|\eta_{0}^r\|_2 \leq c \sqrt{s/n},    \ \|\eta^r_{0}\|_1 \leq c \sqrt{s^2/n},\\
  s = c n^{\frac{1}{2\mathsf{a}}},
\end{array}
\end{align}
where the constant $c$ depends only on $(\sa, \sA)$.

\textsc{Condition AS.2}  \textit{We assume that }
\begin{equation}
s^2 \log(pn)^2/n \leq o(1).
\end{equation}

We shall perform inference on $\alpha_0$ using the empirical analog of theoretical
equations:
\begin{align}\label{moment}
{\M}(\alpha_0, \eta_0) = 0,  \quad {\M}(\alpha, \eta) := \textrm{E} \left[ \psi(w_i, \alpha, \eta)\right],
\end{align}
where $\psi= (\psi_k)_{k=1}^{p^d}$ is defined by
$$\psi_{k}(w_i, \alpha, \eta) := \left (y_i - x_{i}'\theta-  \sum_{\bar k=1}^{p^d}(d_{i \bar k}-x_{i}' \vartheta_{\bar k}) \alpha_{\bar k}\right)
(x_{i}'  \gamma_k+ z_{i}' \delta_k - x_{i}'\vartheta_k).$$

We can verify that the following orthogonality condition holds:
\begin{align}
\label{orthogonality} \partial_{\eta'}  {\M}(\alpha_0, \eta) \Big \vert_{\eta= \eta_0} =0.
\end{align}
This means that missing the true value $\eta_0$ by a small amount does not invalidate the moment condition.
Therefore, the  moment condition will be relatively insensitive to
non-regular estimation of $\eta_0$.

We denote the empirical analog of (\ref{moment}) as
\begin{align}
\label{emp moment}
\hat{\M}(\alpha, \hat \eta) = 0,  \quad  \hat{\M}(\alpha, \eta)  := \En  \left[ \psi_i(\alpha, \eta)\right].
\end{align}
Inference based on this condition can be shown to be immunized against small selection mistakes by virtue of orthogonality.

The above formulation is a special case of the linear-affine model.  Indeed, here we have
$$
\M(\alpha, \eta) =  \Gamma_1(\eta) \alpha + \G_2(\eta), \quad \hat{\M}(\alpha, \eta) = \hat \Gamma_1(\eta) \alpha + \hat \Gamma_2(\eta),
$$ $$
\Gamma_1(\eta) = \Ep[ \psi^a(w_i, \eta )], \quad \hat \Gamma_1(\eta) = \En[\psi^a(w_i, \eta ) ],
$$ $$
\Gamma_2(\eta) = \Ep[ \psi^b(w_i, \eta )], \quad \hat \Gamma_2(\eta) = \En[\psi^b(w_i, \eta ) ],
$$
where $$\psi^a_{k, \bar k}(w_i, \eta) = - (d_{i\bar k}-x_{i}' \vartheta_{\bar k})(x_{i}'  \gamma_k+ z_{i}' \delta_k - x_{i}'\vartheta_k) ,$$$$\psi^b_{k}(w_i, \eta) =  (y_i - x_{i}'\theta)(x_{i}'  \gamma_k+ z_{i}' \delta_k - x_{i}'\vartheta_k).$$



Consequently we can use the results of the previous section. In order to do so we need to
provide a suitable estimator for $\eta_0$.  Here we use the Lasso and Post-Lasso estimators,
as defined in \cite{BCCH12}, to deal with non-normal errors and heteroscedasticity.

\begin{algorithm}[Estimation of $\eta_0$] (1) For each $k$, do Lasso or Post-Lasso Regression of $d_{ik}$ on $x_i, z_i$ to obtain $\hat \gamma_k$ and $\hat \delta_k$. (2) Do Lasso or Post-Lasso Regression of $y_i$ on $x_i$ to get $\hat \theta$. (3) Do Lasso or Post-Lasso Regression of $\hat d_{ik} = x_i'\hat\gamma_k + z_i'\hat\delta_k$ on $x_i$ to get $\hat \vartheta_k$. The estimator of $\eta_0$ is given by $\hat \eta = ({\hat \theta}', ({\hat \vartheta}_{k}', \hat \gamma_{0k}', \hat \delta_{k}')_{k=1}^{p^d} )'.$
 \end{algorithm}

We then use
$$
\hat \Omega(\alpha, \hat \eta) = \En [\psi(w_i, \alpha, \hat \eta) \psi(w_i, \alpha, \hat \eta)'].
$$
to estimate the variance matrix $\Omega(\alpha, \eta_0)= \En [\psi(w_i, \alpha, \eta_0) \psi(w_i, \alpha, \eta_0)'].$
We formulate the orthogonal score statistic and the $C(\alpha)$-statistic,
\begin{equation}\label{eq:normality3}
S(\alpha) :=  \hat \Omega^{-1/2}_n(\alpha, \hat \eta)\sqrt{n}\hat{\M}(\alpha, \hat \eta), \quad
C(\alpha) = \| S(\alpha)\|^2,
\end{equation}
as well as our estimator $\hat \alpha$:
$$
\hat \alpha = \arg\min_{\alpha \in \AA} \| \sqrt{n}\hat{\M}(\alpha, \hat \eta) \|^2.
$$
Note also that $\hat \alpha = \arg\min_{\alpha \in \AA} C(\alpha)$ under mild conditions, since we work
with ``exactly identified" systems of equations.   We also need to specify a variance estimator
$\hat V_n$ for the large sample variance $V_n$ of $ \hat \alpha$.  We set $\hat V_n = (\hat \Gamma_1(\hat \eta)')^{-1}
\hat \Omega(\hat \alpha, \hat \eta) (\hat \Gamma_1(\hat \eta))^{-1}$.

To estimate the nuisance parameter  we impose the following condition.
Let  $f_{i} := (f_{ij})_{j=1}^{p_{f}} := (x_i', z_i')'$; $h_{i} : = (h_{il})_{l=1}^{p_h} := (y_i, d'_i,  \bar d'_i)'$
where $\bar d_i = (\bar d_{ik})_{k=1}^{p^d}$ and $\bar d_{ik} := x_i'\gamma_{0k} + z_i'\delta_{0k}$;
 $v_{i} = (v_{il})_{l=1}^{p_h} := (\varepsilon_i, \rho^y_i, {\rho_i^{d}}',  {\varrho_i}')'$
   where $\varrho_i = (\varrho_{ik})_{k=1}^{p^d}$ and $\varrho_{ik} := d_{ik} - \bar d_{ik}$.
   Let $\tilde h_{i} := h_i - \Ep[h_i]$.

\textsc{Condition RF.}  \textit{  (i) The eigenvalues of $\Ep[f_i f_i']$ are bounded from
above by $\sC$ and from below by $\sc$. For all $j$ and $l$,
(ii) $\Ep[h^2_{il}] +\Ep[ |f^2_{ij} \tilde h^2_{il}|] + 1/\Ep[f_{ij}^2 v_{il}^2]  \leq \sC$ and
$\Ep[ |f^2_{ij} v^2_{il}|]  \leq \Ep[ |f^2_{ij} \tilde h^2_{il}|]$, (iii) $
\Ep[|f_{ij}^3 v_{il}^3|]^2 \log^3 (p n)/n \leq \delta_n$, and (iv) $s \log (pn)/n \leq \delta_n$.
With probability no less than $1- \delta_n$, we have that (v)  $  \max_{i \leq n, j}  f_{ij}^2 [s^2 \log (pn)]/n \leq \delta_n
$ and $ \max_{l, j }   |(\En - \Ep)[f_{ij}^2 v^2_{il}] |  + | (\En - \Ep)[f^2_{ij}\tilde h^2_{il}]| \leq \delta_n$
and (vi)  $\|\En[f_i f_i']- \Ep[f_i f_i']\|_{\mathsf{sp}(\ell_n s)} \leq \delta_n$.}

The conditions are motivated by those given in \cite{BCCH12}.
The current conditions are made slightly stronger to account
for the fact that we use zero covariance conditions in formulating the moments.  Some conditions
could be easily relaxed at a cost of more complicated exposition.

To estimate the variance matrix and establish asymptotic normality, we also need
the following condition. Let $\sq>4$ be a fixed constant.

\textsc{Condition SM.} \textit{  For each $l$ and $k$, (i) $\displaystyle \Ep [|h_{il}|^\sq] +  \Ep[|v_{il}|^\sq] \leq \sC$, (ii) $ \sc \leq \Ep[\varepsilon_i^2\mid x_i, z_i]\leq \sC$,  $\sc< \Ep[{\varrho_{ik}^2}\mid x_i,z_i]\leq \sC$  a.s.,  (iii) $\sup_{\alpha \in \AA} \|\alpha\|_2 \leq \sC$.}

Under the conditions set forth above, we have the following result on validity of post-selection
and post-regularization
inference using the $C(\alpha)$-statistic and estimators derived from it.

\begin{proposition}[Valid  Inference in Large Linear Models using $C(\alpha)$-statistics]
\label{prop: IV Calpha} Let $\mathbf{P}_n$ be the collection of all $\Pr$ such that Conditions AS.1-2, SM, and RF hold for the given $n$.   Then uniformly in $\Pr \in \mathbf{P}_n$,
$S(\alpha_0) \leadsto \mathcal{N}(0, I)$, and $C(\alpha_0) \leadsto \chi^2(p^d)$. As a consequence,
the confidence set $\mathsf{CR}_{1-a} = \{ \alpha \in \mathcal{A}:  C(\alpha) \leq  c(1-a) \}$, where
$c(1-a)$ is the $1-a$-quantile of a $\chi^2(p^d)$  is uniformly valid for $\alpha_0$, in the sense that
$$
\lim_{n \to \infty} \sup_{\Pr \in \mathbf{P}_n}  |\Pr (  \alpha_0 \in \mathsf{CR}_{1-a} ) - (1-a) |  = 0.
$$
Furthermore, for $V_n =  (\Gamma_1')^{-1}   \Omega(\alpha_0,\eta_0) (\Gamma_1)^{-1},$
we have that
$$
\lim_{n \to \infty} \sup_{\Pr \in \mathbf{P}_n} \sup_{R \in \mathcal{R}} |\Pr (  V_n^{-1/2} (\hat \alpha - \alpha_0) \in R) -  \mathbb{P}( \mathcal{N}(0,I) \in R) |  = 0,
$$
where $\mathcal{R}$ is the collection of all convex sets. Moreover, the result continues to apply if $V_n$ is replaced by $\hat V_n$.   Thus,  $ \mathsf{CR}^l_{1-a} = [
l'\hat \alpha \pm c(1-a/2) (l'\hat V_nl/n)^{1/2}]$, where $c(1-a/2)$ is the $(1-a/2)$-quantile
of a $\mathcal{N}(0,1)$, provides a uniformly valid confidence set for $l'\alpha_0$:
$$
\lim_{n \to \infty} \sup_{\Pr \in \mathbf{P}_n}   |\Pr (  l'\alpha_0 \in \mathsf{CR}^l_{1-a} ) - (1-a) |  = 0.
$$

\end{proposition}

\subsection{Simulation Illustration}

In this section, we provide results from a small Monte Carlo simulation to illustrate the performance of the estimator resulting from the application of Algorithm 1 in a small sample setting.  As comparison, we report results from two commonly used
``unprincipled"  alternatives for which uniformly valid inference over the class of approximately sparse models does not hold.  Simulation parameters were chosen so that approximate sparsity holds but exact sparsity is violated in such a way that we expected the unprincipled procedures to perform poorly.

For our simulation, we generate data as $n$ iid draws from the model
\begin{align*}
\left. \begin{array}{ll}
y_i &= \alpha d_i + x_i' \beta + 2\varepsilon_i \\
d_i &= x_{i}'\gamma + z_{i}'\delta + u_i \\
z_i &= \Pi x_i + .125\zeta_i \\
\end{array} \right |  \quad \left( \begin{array}{c} \varepsilon_i \\ u_i \\ \zeta_i \\ x_i \end{array} \right )
\sim \mathcal{N}\left(0 , \left( \begin{array}{cccc} 1 & .6 & 0 & 0 \\ .6 & 1 & 0 & 0 \\ 0 & 0 & I_{p_n^z} & 0 \\ 0 & 0 & 0 & \Sigma \end{array} \right) \right),
\end{align*}
where $\Sigma$ is a $p_n^x \times p_n^x$ matrix with $\Sigma_{kj} = (0.5)^{|j-k|}$ and $I_{p_n^z}$ is a $p_n^z \times p_n^z$ identity matrix.  We set the number of potential controls variables ($p_n^x$) to 200, the number of instruments ($p_n^z$) to 150, and the number of observations ($n$) to 200.  For model coefficients, we set $\alpha = 0$, $\beta = \gamma$ as $p_n^x-$vectors with entries $\beta_j = \gamma_j =  1/(9\nu)$, $\nu = {4/9 + \sum_{j = 5}^{p_n^x} 1/j^2}$ for $j \le 4$ and $\beta_j = \gamma_j = 1/(j^2\nu)$ for $j > 4$, $\delta$ as a $p_n^z-$vector with entries $\delta_j = \frac{3}{j^2}$, and $\Pi = [I_{p_n^z} \ , \ 0_{p_n^z \times (p_n^x-p_n^z)}].$  We report results based on 1000 simulation replications.

We provide results for four different estimators - an infeasible Oracle estimator that knows the nuisance parameters $\eta$ (Oracle), two naive estimators, and the proposed ``Double-Selection'' estimator.  The results for the proposed ``Double-Selection'' procedure are obtained following Algorithm 1 using Post-Lasso at every step.  To obtain the Oracle results, we run standard IV regression of $y_i - \textrm{E}[y_i|x_i]$ on $d_i-\textrm{E}[d_i|x_i]$ using the single instrument $\zeta_i' \delta$.  The expected values are obtained from the model above and $\zeta_i' \delta$ provides the information in the instruments that is unrelated to the controls.  

The two naive alternatives offer unprincipled, although potentially intuitive alternatives.  The first naive estimator follows Algorithm 1 but replaces Lasso/Post-Lasso with stepwise regression with a p-value for entry of .05 and a p-value for removal of .10 (Stepwise).   The second naive estimator (Non-orthogonal) corresponds to using a moment condition that does not satisfy the orthogonality condition described previously but will produce valid inference when perfect model selection in the regression of $d$ on $x$ and $z$ is possible or perfect model selection in the regression of $y$ on $x$ is possible and an instrument is selected in the $d$ on $x$ and $z$ regression.\footnote{Specifically, for the second naive alternative (Non-orthogonal), we first do Lasso regression of $d$ on $x$ and $z$ to obtain Lasso estimates of the coefficients $\gamma$ and $\delta$.  Denote these estimates as $\hat\gamma_L$ and $\hat\delta_L$, and denote the indices of the coefficients estimated to be non-zero as $\hat{I}^d_x = \{j : \hat\gamma_{Lj} \ne 0\}$ and $\hat{I}^d_z = \{j : \hat\delta_{Lj} \ne 0\}$. We then run Lasso regression of $y$ on $x$ to learn the identities of controls that predict the outcome. We denote the Lasso estimates as $\hat\theta_L$ and keep track of the indices of the coefficients estimated to be non-zero as $\hat{I}^y_x = \{j : \hat\theta_{Lj} \ne 0\}$. We then take the union of the controls selected in either step $\hat{I}_x = \hat{I}^y_x \cup \hat{I}^d_x$.  The estimator of $\alpha$ is then obtained as the usual 2SLS estimator of $y_i$ on $d_i$ using all selected elements from $x_i$, $x_{ij}$ such that $j \in \hat{I}_x$, as controls and the selected elements from $z_i$, $z_{ij}$ such that $j \in \hat{I}^d_z$, as instruments. }

All of the Lasso and Post-Lasso estimates are obtained using the data-dependent penalty level from \citen{BC-PostLASSO}.  This penalty level depends on a standard deviation that is estimated adapting the iterative algorithm described in \citen{BCCH12} Appendix A using Post-Lasso at each iteration.  For inference in all cases, we use standard t-tests based on conventional homoscedastic IV standard errors obtained from the final IV step performed in each strategy.

\begin{figure}\label{SimulationDistributions}
\includegraphics[width=\columnwidth]{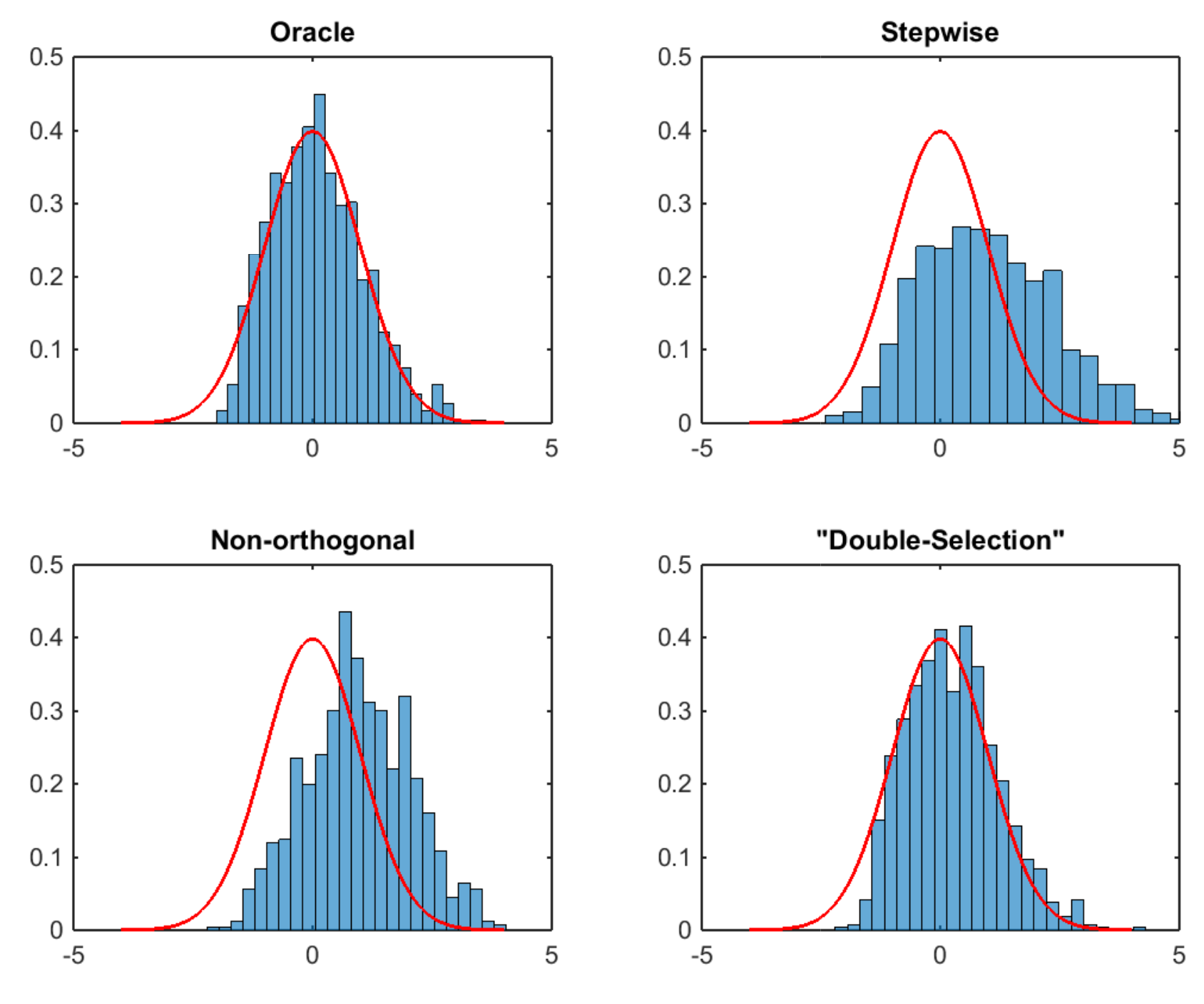}
\caption{The figure presents the histogram of the estimator from each method centered around the true parameters and scaled by the estimated standard error from the simulation experiment.  The red curve is the pdf of a standard normal which will correspond to the sampling distribution of the estimator under the asymptotic approximation.  Each panel is labeled with the corresponding estimator from the simulation.}
\end{figure}

We display the simulation results in Figure \ref{SimulationDistributions}, and we report the median bias (Bias), median absolute deviation (MAD), and size of 5\% level tests (Size) for each procedure in Table \ref{SimulationTable}.  For each estimator, we plot the simulation estimate of the sampling distribution of the estimator centered around the true parameter and scaled by the estimated standard error.  With this standardization, usual asymptotic approximations would suggest that these curves should line up with a $\mathcal{N}(0,1)$ density function which is displayed as the bold solid line in the figure.  We can see that the Oracle estimator and the Double-Selection estimator are centered correctly and line up reasonably well with the $\mathcal{N}(0,1)$, although both estimators exhibit some mild skewness.  It is interesting that the sampling distributions of the Oracle and Double-Selection estimators are very similar as predicted by the theory.  In contrast, both of the naive estimators are centered far from zero, and it is clear that the asymptotic approximation provides a very poor guide to the finite sample distribution of these estimators in the design considered.

\begin{table}
\begin{center}
\caption{Summary of Simulation Results for the Estimation of $\alpha$}\label{SimulationTable}
\begin{tabular}{l c c c}
\hline \hline
Method & Bias &  MAD     & Size   \\ 
\hline
Oracle	&	0.015	&	0.247	&	0.043	\\
Stepwise	&	0.282	&	0.368	&	0.261	\\
Non-orthogonal	&	0.084	&	0.112	&	0.189 \\
Double-Selection	&	0.069	&	0.243	&	0.053	\\
\hline
\hline
\end{tabular}
\end{center}
\begin{flushleft}
\footnotesize{This table summarizes the simulation results from a linear IV model with many instruments and controls.  Estimators include an infeasible oracle as a benchmark (Oracle), two naive alternatives (Stepwise and Non-orthogonal) described in the text, and our proposed feasible valid procedure (Double-Selection).  Median bias (Bias), median absolute deviation (MAD), and size for 5\% level tests (Size) are reported.}
\end{flushleft}
\end{table}

The poor inferential performance of the two naive estimators is driven by different phenomena.  The unprincipled use of stepwise regression fails to control spurious inclusion of irrelevant variables which leads to inclusion of many essentially irrelevant variables, resulting in many-instrument-type problems (e.g. \citen{NeweyEtAl-JIVE}).  In addition, the spuriously included variables are those most highly correlated to the noise within sample which adds an additional type of ``endogeneity bias''.  The failure of the ``Non-orthogonal'' method is driven by the fact that perfect model selection is not possible within the present design: Here we have
model selection mistakes in which the control variables that are correlated to the instruments but only moderately correlated to the outcome and endogenous variable are missed.  Such exclusions result in standard omitted variables bias in the estimator for the parameter of interest and substantial size distortions.  The additional step in the Double-Selection procedure can be viewed as a way to guard against such mistakes.   Overall, the results illustrate the uniformity claims made in the preceding section.  The feasible Double-Selection procedure following from Algorithm 1 performs similarly to the semi-parametrically efficient infeasible Oracle.  We obtain good inferential properties with the asymptotic approximation providing a fairly good guide to the behavior of the estimator despite working in a setting in which perfect model selection is impossible.  Although simply illustrative of the theory, the results are reassuring and in line with extensive simulations in the linear model with many controls provided in \citen{BelloniChernozhukovHansen2011}, in the instrumental variables model with many instruments and a small number of controls provided in \citen{BCCH12}, and in linear panel data models provided in \citen{BCHK:FE}.

\subsection{Empirical Illustration: Logit Demand Estimation}

As further illustration of the approach, we provide a brief empirical example in which we estimate the coefficients in a simple logit model of demand for automobiles using market share data.  Our example is based on the data and most basic strategy from \citeasnoun{BLP}.  Specifically, we estimate the parameters from the model 
\begin{align*}
\log(s_{it}) - \log(s_{0t}) &= \alpha_0 p_{it} + x_{it}'\beta_0 + \varepsilon_{it}, \\
p_{it} &= z_{it}'\delta_0 + x_{it}'\gamma_0 + u_{it},
\end{align*}
where $s_{it}$ is the market share of product $i$ in market $t$ with product zero denoting the outside option, $p_{it}$ is price and is treated as endogenous, $x_{it}$ are observed included product characteristics, and $z_{it}$ are instruments.  One could also adapt the proposed variable selection procedures to extensions of this model such as the nested logit model or models allowing for random coefficients; see, e.g., \citeasnoun{BLPLASSO} for an example with a random coefficient.  

In our example, we use the same set of product characteristics ($x$-variables) as used in obtaining the basic results in \citen{BLP}.  Specifically, we use five variables in $x_{it}$: a constant, an air conditioning dummy, horsepower divided by weight, miles per dollar, and vehicle size.  We refer to these five variables as the baseline set of controls.

We also adopt the argument from \citen{BLP} to form our potential instruments.  \citen{BLP} argue that that characteristics of other products will satisfy an exclusion restriction, $\textrm{E}[\varepsilon_{it}|x_{j\tau}] = 0$ for any $\tau$ and $j \ne i$, and thus that any function of characteristics of other products may be used as instrument for price.  This condition leaves a very high-dimensional set of potential instruments as any combination of functions of $\{x_{j\tau}\}_{j \ne i, \tau \ge 1}$ may be used to instrument for $p_{it}$.  To reduce the dimensionality, \citeasnoun{BLP} use intuition and an exchangeability argument to motivate consideration of a small number of these potential instruments formed by taking sums of product characteristics formed by summing over products excluding product $i$.  Specifically, we form baseline instruments by taking
$$
z_{k,it} = \left(\sum_{r \ne i, r \in \mathcal{I}_f} x_{k,rt} , \sum_{r \ne i, r \notin \mathcal{I}_f} x_{k,rt} \right) 
$$
where $x_{k,it}$ is the $k^{th}$ element of vector $x_{it}$ and $\mathcal{I}_f$ denotes the set of products produced by firm $f$. 
This choice yields a vector $z_{it}$ consisting of 10 instruments.  We refer to this set of instruments as the baseline instruments.

Although the choice of the baseline instruments and controls is motivated by good intuition and economic theory, we note that theory does not clearly state which product characteristics or instruments should be used in the model.  Theory also fails to indicate the functional form with which any such variables should enter the model.  The high-dimensional methods outlined in this paper offer one strategy to help address these concerns that complements the economic intuition motivating the baseline controls and instruments.  As an illustration, we consider an expanded set of controls and instruments.  We augment the set of potential controls with all first order interactions of the baseline variables, quadratics and cubics in all continuous baseline variables, and a time trend that yields a total of 24 $x$-variables.  We refer to these as the augmented controls.  We then take sums of these characteristics as potential instruments following the original strategy which yields 48 potential instruments.  

We present estimation results in Table \ref{EmpiricalTable}.  We report results obtained by applying the method outlined in Algorithm 1 using just the baseline set of five product characterstics and 10 instruments in the row labeled ``Baseline 2SLS with Selection'' and results obtained by applying the method to the augmented set of 24 controls and 48 instruments in the row labeled ``Augmented 2SLS with Selection.''  In each case, we apply the method outlined in Algorithm 1 using post-Lasso in each step and forcing the intercept to be included in all models.  We employ the heteroscedasticity robust version of Post-Lasso of \citen{BCCH12} following the implementation algorithm provided in Appendix A of \citen{BCCH12}.  For comparison, we also report OLS and 2SLS estimates using only the baseline variables in ``Baseline OLS'' and ``Baseline 2SLS,'' respectively; and we report  OLS and 2SLS estimates using the augmented variable set in ``Augmented OLS'' and ``Augmented 2SLS,'' respectively.   All standard errors are conventional heteroscedasticity robust standard errors. 

\begin{table}
\begin{center}
\caption{Estimates of Price Coefficient}\label{EmpiricalTable}
\begin{tabular}{l c c c}
\hline\hline
 & Price Coefficient &  Standard Error  &  Number Inelastic    \\
\cline{2-4}
\text{} & \multicolumn{3}{c}{\textit{Estimates Without Selection}} \\
Baseline OLS	&	-0.089	&	0.004	& 1502   \\
Baseline 2SLS	&	-0.142	&	0.012	&	670	\\
Augmented OLS	&	-0.099	&	0.005	&	 1405 \\
Augmented 2SLS	&	-0.127	&	0.014	&	874	\\
\text{} & \multicolumn{3}{c}{\textit{2SLS Estimates With ``Double Selection"}} \\
Baseline 2SLS  Selection & -0.185  & 0.014      & 139    \\
Augmented 2SLS Selection &  -0.221    & 0.015      & 12     \\
\hline
\hline
\end{tabular}
\end{center}
\begin{flushleft}
\footnotesize{This table reports estimates of the coefficient on price (``Price Coefficient'') along with the estimated standard error (``Standard Error'') obtained using different sets of controls and instruments.  The rows ``Baseline OLS'' and ``Baseline 2SLS'' respectively provide OLS and 2SLS results using the baseline set of variables (5 controls and 10 instruments) described in the text.  The rows ``Augmented OLS,'' ``Augemented 2SLS ''are defined similarly but use the augmented set of variables described in the text (24 controls and 48 instruments).    The rows ``Baseline 2SLS with Selection'' and ``Augmented 2SLS with Selection'' 
 applies the ``double selection" approach developed in this paper to select a set of controls and instruments and perform valid post-selection inference about the estimated price coefficient where selection occurs considering only the baseline variables.   For each procedure, we also report the point estimate of the number of products for which demand is estimated to be inelastic in the column ``Number Inelastic.''}
\end{flushleft}
\end{table}

Considering first estimates of the price coefficient, we see that the estimated price coefficient increases in magnitude as we move from OLS to 2SLS and then to the selection based results.  After selection using only the original variables, we estimate the price coefficient to be -.185 with an estimated standard error of .014 compared to an OLS estimate of -.089 with estimated standard error of .004 and 2SLS estimate of -.142 with estimated standard error of .012.  In this case, all five controls are selected in the log-share on controls regression, all five controls but only four instruments are selected in the price on controls and instruments regression, and four of the controls are selected for the price on controls relationship.  The difference between the baseline results is thus largely driven by the difference in instrument sets.  The change in the estimated coefficient is consistent with the wisdom from the many-instrument literature that inclusion of irrelevant instruments biases 2SLS toward OLS.  

With the larger set of variables, our post-model-selection estimator of the price coefficient is -.221 with an estimated standard error of .015 compared to the OLS estimate of -.099 with an estimated standard error of .005 and 2SLS estimate of -.127 with an estimated standard error of .014.  Here, we see some evidence that the original set of controls may have been overly parsimonious as we select some terms that were not included in the baseline variable set.  We also see a closer agreement between the OLS estimate and 2SLS estimate without selection which is likely driven by the larger number of instruments considered and the usual bias towards OLS seen in 2SLS with many weak or irrelevant instruments.  In the log-share on controls regression, we have eight control variables selected; and we have seven controls and only four instruments selected in the price on controls and instrument regression.  We also have 13 variables selected for the price on controls relationship.  The selection of these additional variables suggests that there is important nonlinearity missed by the baseline set of variables.

The most interesting feature of the results is that estimates of own-price elasticities become more plausible as we move from the baseline results to the results based on variable selection with a large number of controls.  Recall that facing inelastic demand is inconsistent with profit maximizing price choice within the present context, so theory would predict that demand should be elastic for all products.  However, the baseline point estimates imply inelastic demand for 670 products.  When we use the larger set of instruments without selection, the number of products for which we estimate inelastic demand increases to 874 with the increase generated by the 2SLS coefficient estimate moving back towards the OLS estimate. 
The use of the variable selection results provides results closer to the theoretical prediction.  The point estimates based on selection from only the baseline variables imply inelastic demand for 139 products, and we estimate inelastic demand for only 12 products using the results based on selection from the larger set of variables.  Thus,
the new methods provide the most reasonable estimates of own-price elasticities.

We conclude by noting that the simple specification above suffers from the usual drawbacks of the logit demand model.  However, the example illustrates how the application of the methods we outlined may be used in the estimation of structural parameters in economics and adds to the plausibility of the resulting estimates.  In this example, we see that we obtain more sensible estimates of key parameters with at most a modest cost in increased estimation uncertainty after applying the methods in this paper while considering a flexible set of variables.

\section{Overview of Related Literature}

Inference following model selection or regularization more generally has been an active area of research in econometrics and statistics for the last several years.  In this section, we provide a brief overview of this literature highlighting some key developments.  This review is necessarily selective due to the large number of papers available and the rapid pace at which new papers are appearing.  We choose to focus on papers that deal specifically with high-dimensional nuisance parameter settings, and note that the ideas in these papers apply in low dimensional settings as well.

Early work on inference in high-dimensional settings focused on inference based on the so-called perfect recovery; see, e.g., \citen{FanLi2001} for an early paper, \citen{FanLv2010:Review} for a more recent review, and \citen{bulman:book} for a textbook treatment. A consequence of this property is that model selection does not impact the asymptotic distribution of the parameters estimated in the selected model.  This feature allows one to do inference using standard approximate distributions for the parameters of the selected model ignoring that model selection was done.  While convenient and fruitful in many applications (e.g. signal processing), such results effectively rely on strong conditions that imply that one will be able to perfectly select the correct model.  For example, such results in linear models require the so called ``beta-min condition" (\cite{bulman:book}) that all but a small number of coefficients are exactly zero and the remaining non-zero coefficients are bounded away from zero, effectively ruling out variables that have small, non-zero coefficients.  Such conditions seem implausible in many applications, especially in econometrics, and relying on such conditions produces asymptotic approximations that may provide very poor approximations to finite-sample distributions of estimators as they are not uniformly valid over sequences of models that include even minor deviations from conditions implying perfect model selection.  The concern about the lack of uniform validity of inference based on oracle properties was raised in a series of papers, including \citen{leeb:potscher:review} and \citen{leeb:potscher:hodges} among many others,
and the more recent work on post-model-selection inference has been focused on offering procedures that provide uniformly valid inference over interesting (large) classes of models that include cases where perfect model selection will not be possible.

To our knowledge, the first work to formally and expressly address the problem of obtaining uniformly valid inference following model selection is \citen{BellChernHans:Gauss} which considered inference about parameters on a low-dimensional set of endogenous variables following selection of instruments from among a high-dimensional set of potential instruments in a homoscedastic, Gaussian instrumental variables (IV)  model.  The approach does not rely on implausible ``beta-min" conditions which imply perfect model selection but instead relies on the fact that the moment condition underlying IV estimation satisfies the \textit{orthogonality condition} (\ref{eq:ortho}) and the use of high-quality variable selection methods.  These ideas were further developed in the context of providing uniformly valid inference about the parameters on endogenous variables in the IV context with many instruments to allow non-Gaussian heteroscedastic disturbances in \citen{BCCH12}.  These principles have also been applied in \citen{BCH2011:InferenceGauss}, who developed approaches for regression and IV models with Gaussian errors; \citen{BelloniChernozhukovHansen2011} (ArXiv 2011), which covers estimation of the parametric components of the partially linear model, estimation of average treatment effects, and provides a formal statement of the orthogonality condition (\ref{eq:ortho}); \citen{Farrell:JMP} which covers average treatment effects with discrete, multi-valued treatments; \citen{Kozbur:JMP} which covers additive nonparametric models; and \citen{BCHK:FE} which extends the IV and partially linear model results to allow for fixed effects panel data and clustered dependence structures.  The most recent, general approach is provided in \citen{BCFH:Policy} where inference about parameters defined by a continuum of orthogonalized estimating equations with infinite-dimensional nusiance parameters is analyzed and positive results on inference are developed.  The framework in \citen{BCFH:Policy} is general enough to cover the aforementioned papers and many other parametric and semi-parametric models considered in economics.

As noted above, providing uniformly valid inference following model selection is closely related to use of Neyman's $C(\alpha)$-statistic. Valid confidence regions can be obtained by inverting tests based on these statistics, and minimizers of $C(\alpha)$-statistics may be used as point estimators.  The use of $C(\alpha)$ statistics for testing and estimation in high-dimensional approximately sparse models was first explored in the context of high-dimensional quantile regression in \citen{BCK-LAD} (Oberwolfach, 2012) and \citen{BCK-QR} and in the context of high-dimensional logistic regression and other high-dimensional generalized linear models by \citen{BCY-honest}.  More recent uses of $C(\alpha)$-statistics (or close variants, under different names) include those in \citen{witten:score}, \citen{hanliu1}, and \citen{hanliu2} among others.

There have also been parallel developments based upon ex-post ``de-biasing'' of estimators.  This approach is mathematically equivalent to doing classical ``one-step" corrections in the general framework of Section 2.
 Indeed, while at first glance this ``de-biasing" approach may appear distinct from that taken in the papers listed above in this section, it is the same as approximately solving -- by doing one Gauss-Newton step -- orthogonal estimating equations satisfying (\ref{eq:ortho}).
 The general results of Section 2 suggest that these approaches -- the exact solving and ``one-step" solving -- are generally first-order asymptotically equivalent, though higher-order differences may persist. To the best of our knowledge, the ``one-step" correction approach was first
 employed in high-dimensional sparse models by \citen{ZhangZhang:CI}  (ArXiv 2011) which covers the homoscedastic linear model (as well as in several follow-up works by the authors).    This approach has been further used in \citen{vdGBRD:AsymptoticConfidenceSets} (ArXiv 2013) which covers homoscedastic linear models and some generalized linear models, and \citen{JM:ConfidenceIntervals} (ArXiv 2013) which offers a related, though somewhat different approach.
 Note that \citen{BCK-LAD} and \citen{BCY-honest} also offer results on ``one-step" corrections as part of their analysis
 of estimation and inference based upon the orthogonal estimating equations.  We would not expect that the use of orthogonal estimating equations or the use of ``one-step" corrections to dominate each other in all cases, though computational evidence in \citen{BCY-honest} suggests that the use of exact solutions to orthogonal estimating equations may be preferable to approximate solutions obtained from ``one-step" corrections in the contexts considered in that paper.

Another branch of the recent literature takes a complementary, but logically distinct, approach that aims at doing valid inference for the parameters of a ``pseudo-true'' model that results from the use of a model selection procedure, see \citen{BerkEtAl13}.   
Specifically, this approach conditions on a model selected by a data-dependent rule and then attempts to do  inference --  conditional on the selection event -- for the parameters of the selected model, which may deviate from the ``true'' model that generated the data.   Related developments within this approach appear in \citen{AdaptiveGraphLasso}, \citen{LeeTaylorScreening}, \citen{LeeEtAlLasso}, \citen{LockhartEtAlLasso}, \citen{LoftusStepwise}, \citen{TaylorEtAlAdaptive}, and \citen{FST14}. 
It seems intellectually very interesting to combine the developments of the present paper (and other preceding papers cited above) with developments in this literature.

The previously mentioned work focuses on doing inference for low dimensional parameters in the presence of high dimensional nuisance parameters.  There have also been developments on performing inference for high dimensional parameters.
\cite{Chern:SG} proposed inverting a Lasso performance bound in order to construct a simultaneous, Scheff\'{e}-style
confidence band on all parameters. An interesting feature of this approach is that it uses
weaker design conditions than many other approaches but requires the data analyst to supply explicit bounds on restricted eigenvalues.  \cite{GautierTsybakovHDIV} (ArXiv 2011)
and \cite{CCK:AOS13} employ similar ideas while also working with various generalizations of restricted eigenvalues.
\cite{vdg:nickl} construct confidence ellipsoids for the entire parameter vector using sample splitting ideas.
Somewhat related to this literature are the results of \citen{BCK-LAD} who use the orthogonal
estimating equations framework with infinite-dimensional nuisance parameters and construct a simultaneous confidence rectangle
for many target parameters where the number of target parameters could be much larger than the sample size.
They relied upon the high-dimensional central limit theorems and bootstrap results established in \citen{CCK:AOS13}.

Most of the aforementioned results rely on (approximate) sparsity and related sparsity-based estimators.  Some examples of the use of alternative regularization schemes are available in the many instrument literature in econometrics.  For example, \citen{chamberlain:imbens:reiv} use a shrinkage estimator resulting from use of a Gaussian random coefficients structure over first-stage coefficients, and \citen{okui:manyiv} uses ridge regression for estimating the first-stage regression in a framework where the instruments may be ordered in terms of relevance. \citen{carrasco:regularizedIV} employs a different strategy based on directly regularizing the inverse that appears in the definition of the 2SLS estimator allowing for a number of moment conditions that are larger than the sample size; see also \citen{carrasco:LIML}.  The theoretical development in \citen{carrasco:regularizedIV} relies on restrictions on the covariance structure of the instruments rather than on the coefficients of the instruments.  \citen{RJIVE} considers a combination of ridge-regularization and the jackknife to provide a procedure that is valid allowing for the number of instruments to be greater than the sample size under weak restrictions on the covariance structure of the instruments and the first-stage coefficients. In all cases, the orthogonality condition holds allowing root-$n$ consistent and asymptotically normal estimation of the main parameter $\alpha$.

Many other interesting procedures beyond those mentioned in this review have been developed for estimating high-dimensional models; see, e.g. \citen{elements:book} for a textbook review.  Developing new techniques for estimation in high-dimensional settings is also still an active area of research, so the list of methods available to researchers continues to expand.  The use of these procedures and the impact of their use on inference about low-dimensional parameters of interest is an interesting research direction to explore.  It seems likely that many of these procedures will provide sufficiently high-quality estimates that they may be used for estimating the high-dimensional nuisance parameters $\eta$ in the present setting.

\appendix

\section{The Lasso and Post-Lasso Estimators in the Linear Model}

Suppose we have data $\{y_i,x_i\}$ for individuals $i = 1,...,n$ where $x_i$ is a $p$-vector of predictor variables and $y_i$ is an outcome of interest.  Suppose that we are interested in a linear prediction model for $y_i$, $y_i = x_i'\eta + \varepsilon_i$, and define the usual least squares criterion function:
$$
\hat Q(\eta):= \frac{1}{n}\sum_{i=1}^{n}(y_{i} - x_i'\eta)^2.
$$
The Lasso estimator is defined as a solution of the following optimization program:
\begin{equation}\label{Def:LASSOmain}
\widehat \eta_{L} \in \arg \min_{\eta \in \mathbb{R}^p} \hat Q (\eta) + \frac{\lambda}{n}\sum_{j=1}^{p} |\psi_j\eta_j|
\end{equation}
where  $\lambda$ is the penalty level and $\{\psi_j\}_{j=1}^{p}$ are covariate specific penalty loadings.  The covariate specific penalty loadings are used to accommodate data that may be non-Gaussian, heteroscedastic, and/or dependent and also help ensure basic equivariance of coefficient estimates to rescaling of the covariates.

The Post-Lasso estimator is defined as the ordinary least square regression applied to the model $\widehat I$ selected by Lasso:\footnote{We note that we can also allow the set $\widehat I$ to contain additional variables not selected by Lasso, but we do not consider that here.}
$$\widehat I = \textrm{support}( \hat \eta_{L} ) = \{ j \in \{1,\ldots,p\} \ : \ |\hat\eta_{L j }| > 0\}.$$
 The Post-Lasso estimator $\hat\eta_{PL}$  is then
\begin{equation}\label{Def:TwoStep} 
\widehat \eta_{PL} \in \arg\min  \{ \hat Q(\eta): \quad \eta \in \mathbb{R}^p \text{ such that } 
\eta_j = 0 \ \text{ for all } j  \notin \widehat I \}
\end{equation}
In words, this estimator is ordinary least squares (OLS) using only the regressors whose coefficients were estimated to be non-zero by Lasso.

Lasso and Post-Lasso are motivated by the desire to predict the target function well without overfitting.  The Lasso estimator is a computationally attractive alternative to some other classic approaches, such as model selection based on information criteria, because it minimizes a convex function.  Moreover, under suitable conditions, the Lasso estimator achieves near-optimal rates in estimating the regression function  $x_i'\eta$.  However, Lasso does suffer from the drawback that the regularization by the $\ell_1$-norm employed in (\ref{Def:LASSOmain}) naturally shrinks all estimated coefficients towards zero causing a potentially significant shrinkage bias. The Post-Lasso estimator is meant to remove some of this shrinkage bias and achieves the same rate of convergence as Lasso under sensible conditions.

Practical implementation of the Lasso requires setting the penalty parameter and loadings.  Verifying good properties of the Lasso typically relies on having these parameters set so that the penalty dominates the score in the sense that
\begin{align*}
\frac{\psi_j\lambda}{n} \geq \max_{j \leq p} 2c \left |\frac{1}{n}\sum_{i=1}^{n} x_{j,i} \varepsilon_i \right | \textrm{ or, equivalently } \frac{\lambda}{\sqrt{n}} \geq \max_{j \leq p} 2c \left |\frac{\frac{1}{\sqrt{n}}\sum_{i=1}^{n} x_{j,i} \varepsilon_i}{\psi_j} \right |
\end{align*}
for some $c > 1$ with high probability.  Heuristically, we would have the term inside the absolute values behaving approximately like a standard normal random variable if we set $\psi_j = \textrm{Var}\left[\frac{1}{\sqrt{n}}\sum_{i=1}^{n} x_{j,i} \varepsilon_i\right]$.  We could then get the desired domination by setting $\frac{\lambda}{2c\sqrt{n}}$ large enough to dominate the maximum of $p$ standard normal random variables with high probability, for example, by setting $\lambda = 2c\sqrt{n}\Phi^{-1}\left(1-{.1/[2p\log(n)]}{}\right)$ where $\Phi^{-1}(\cdot)$ denotes the inverse of the standard normal cumulative distribution function.  Verifying that this heuristic argument holds with large $p$ and data which may not be i.i.d. Gaussian requires careful and delicate arguments as in, for example, \citen{BCCH12} which covers heteroscedastic non-Gaussian data or \citen{BCHK:FE} which covers panel data with within individual dependence.  The choice of the penalty parameter $\lambda$ can also be refined as in \citen{BCW-SqLASSO}.  Finally, feasible implementation requires that $\psi_j$ be estimated which can be done through the iterative procedures suggested in \citen{BCCH12} or \citen{BCHK:FE}.

\section{Proofs}

\subsection{Proof of Proposition \ref{prop:estimation}}  Consider any sequence $\{\Pr_n\}$ in $\{\mathbf{P}_n\}$.

Step 1 ($r_n$-rate).  Here we show that $\| \hat \alpha - \alpha_0 \| \leq r_n$ wp $\to 1$. We have by the identifiability condition, in particular the assumption $\mathrm{mineig}(\Gamma_1'\Gamma_1) \geq c$, that
$$\Pr_n(\|\hat \alpha - \alpha_0\|> r_n) \leq \Pr_n ( \| \M(\hat \alpha, \eta_0)\| \geq \iota(r_n)), \quad
\iota(r_n):= 2^{-1} (\{\sqrt{c} r_n\} \wedge c).$$ Hence it suffices to show that
wp $\to 1$, 
$\| \M(\hat \alpha, \eta_0)\| < \iota(r_n).
$
By the triangle inequality, $$\|\M(\hat \alpha, \eta_0)\| \leq I_1 + I_2 + I_3, \quad
\begin{array}{ll} I_1= \| \M(\hat \alpha, \eta_0) - \M(\hat \alpha, \hat \eta)\|, \\ I_2= \| \M(\hat \alpha, \hat \eta) - \Mn(\hat \alpha, \hat \eta)\|, \\
I_3 =  \|\Mn(\hat \alpha, \hat \eta) \|.
\end{array}$$
By assumption (\ref{eq:se1}), wp $\to 1$
$$
I_1 + I_2 \leq  o(1) \{ r_n  +  I_3 + \|\M(\hat \alpha, \eta_0)\| \}.
$$
Hence,
$$
\|\M(\hat \alpha, \eta_0)\|  (1- o(1))  \leq   o(1)( r_n +    I_3) + I_3.
$$
By construction of the estimator, $$I_3  \leq o(n^{-1/2}) + \inf_{\alpha \in \AA} \| \Mn(\alpha, \hat \eta) \| \lesssim_{\Pr_n} n^{-1/2},$$
which follows because
\begin{equation}\label{eq:minimized objective}
\inf_{\alpha \in \AA} \| \Mn (\alpha, \hat \eta) \| \leq \| \Mn (\bar \alpha, \hat \eta) \| \lesssim_{\Pr_n} n^{-1/2},
 \end{equation}
where $\bar \alpha$ is the one-step estimator defined in Step 3,  as shown in
(\ref{eq:star}).  Hence wp $\to$ 1 $$ \|\M(\hat \alpha, \eta_0)\|   \leq  o(r_n)< \iota(r_n),$$
 where to obtain the last inequality we have used the assumption $\mathrm{mineig}(\Gamma_1'\Gamma_1) \geq c$.

 Step 2 ($n^{-1/2}$-rate).   Here we show that $\| \hat \alpha - \alpha_0 \| \lesssim_{\Pr_n} n^{-1/2}$. By
 condition (\ref{eq:derivatives}) and the triangle inequality,  wp $\to 1$
$$
\| \M(\hat \alpha, \eta_0)\| \geq \| \Gamma_1 (\hat \alpha - \alpha_0)\| - o(1) \| \hat \alpha - \alpha_0\| \geq (\sqrt{c} -o(1)) \| (\hat \alpha - \alpha_0)\| \geq \sqrt{c}/2  \| (\hat \alpha - \alpha_0)\|.
$$
Therefore, it suffices to show that $\| \M(\hat \alpha, \eta_0)\| \lesssim_{\Pr_n} n^{-1/2}$.
We have that  $$\| \M(\hat \alpha, \eta_0)\| \leq II_1+ II_2 + II_3, \quad \begin{array}{lll}
II_1 = \|\M(\hat \alpha, \eta_0) - \M(\hat \alpha, \hat \eta)\|, \\
 II_2=  \|\M(\hat \alpha, \hat \eta) -\Mn(\hat \alpha , \hat \eta)-\Mn(\alpha_0, \eta_0)\|, \\
  II_3 = \| \Mn(\hat \alpha , \hat \eta)\| +  \|\Mn(\alpha_0, \eta_0)\|.
 \end{array} $$
Then, by the orthogonality $\partial_{\eta'} \M(\alpha_0, \eta_0) =0$ and condition (\ref{eq:derivatives}), wp $\to 1$,
\begin{eqnarray*}
II_1 & \leq &  \|\M(\hat \alpha ,\hat \eta) -\M(\hat \alpha , \eta_0) - \partial_{\eta'} \M(\hat \alpha, \eta_0)[\hat \eta - \eta_0]\| +\|\partial_{\eta'} \M(\hat \alpha, \eta_0)[\hat \eta - \eta_0]\| \\
& \leq &  o(1) n^{-1/2} + o(1) \|\hat \alpha - \alpha_0\|  \\
& \leq &  o(1) n^{-1/2} + o(1) (2/\sqrt{c}) \|\M(\hat \alpha, \eta_0)\|.
\end{eqnarray*}
Then, by condition (\ref{eq:se2}) and by $I_3 \lesssim_{\Pr_n} n^{-1/2}$,
\begin{eqnarray*}
II_2 & \leq\ \ \ & o(1)\{n^{-1/2} + \|\Mn(\hat \alpha, \hat \eta)\| + \|\M(\hat \alpha, \eta_0) \|\}\\
 & \lesssim_{\Pr_n}  &   o(1)\{n^{-1/2} + n^{-1/2} + \|\M(\hat \alpha, \eta_0) \|\}.
\end{eqnarray*}
Since $II_3 \lesssim_{\Pr_n} n^{-1/2}$ by (\ref{eq:minimized objective}) and $\|\hat \M(\alpha_0, \eta_0) \| \lesssim_{\Pr_n} n^{-1/2}$ , it follows that wp $\to 1$,
$
 (1-o(1)) \|\M(\hat \alpha, \eta_0) \|   \lesssim_{\Pr_n} n^{-1/2}.
$

Step 3 (Linearization). Define the linearization map $ \alpha \mapsto \Ln(\alpha)$ by
$
\Ln(\alpha):= \Mn(\alpha_0, \eta_0) + \Gamma_1 ( \alpha - \alpha_0).
$
Then $$\|  \Mn(\hat \alpha, \hat \eta) -  \Ln(\hat \alpha) \| \leq III_1 + III_2 +III_3,   
\quad \begin{array}{lll} III_1 = \| \M(\hat \alpha, \hat \eta) - \M(\hat \alpha, \eta_0)  \|,  \\ 
III_2=  \|\M(\hat \alpha, \eta_0) - \Gamma_1(\hat \alpha - \alpha_0)\| ,  \\
  III_3= \|\Mn(\hat \alpha, \hat \eta) - \M(\hat \alpha, \hat \eta) - \Mn( \alpha_0,  \eta_0)\|.
\end{array}
$$
Then, using the assumptions (\ref{eq:derivatives}) and (\ref{eq:se2}), conclude
\begin{eqnarray*}
III_1 & \leq &  \| \M(\hat \alpha, \hat \eta) - \M(\hat \alpha, \eta_0)- \partial_{\eta'} \M (\hat \alpha, \eta_0) [\hat \eta - \eta_0] \|  + \|\partial_{\eta'} \M (\hat \alpha, \eta_0) [\hat \eta - \eta_0]\| \\
& \leq & o(1) n^{-1/2}  +   o(1) \| \hat \alpha - \alpha_0\|, \\
III_2 & \leq & o(1) \| \hat \alpha - \alpha_0\|, \\
III_3 & \leq &  o(1) (n^{-1/2} + \|\Mn(\hat \alpha,\hat \eta)\| + \|\M(\hat \alpha, \eta_0)\|) \\
& \leq &  o(1) (n^{-1/2} +n^{-1/2} + III_2 +  \|\Gamma_1( \hat \alpha - \alpha_0)\|).
\end{eqnarray*}
Conclude that wp $\to 1$, since $\| \Gamma_1'\Gamma_1\| \lesssim 1$ by assumption (\ref{eq:stability}),
$$
\|  \Mn(\hat \alpha, \hat \eta) -  \Ln(\hat \alpha) \|  \lesssim_{\Pr_n}   o(1) (n^{-1/2} + \|\hat \alpha - \alpha_0\|) = o(n^{-1/2}).
$$
Also consider the minimizer of the map $\alpha \mapsto \| \Ln ( \alpha)\|$,  namely,
$$
 \bar \alpha = \alpha_0  - (\Gamma_1' \Gamma_1)^{-1} \Gamma'_1 \hat{\M}(\alpha_0,\eta_0)
$$
which obeys
$
\|\sqrt{n} ( \bar \alpha - \alpha_0)\| \lesssim_{\Pr_n} n^{-1/2}
$
under the conditions of the proposition. We can repeat the argument above to conclude that wp $\to 1$,
$
\|  \Mn(\bar \alpha, \hat \eta) -  \Ln(\bar \alpha) \|  \lesssim_{\Pr_n}   o(n^{-1/2}).
$
This implies, since $\|\Ln(\bar \alpha)\| \lesssim_{\Pr_n} n^{-1/2}$,
\begin{equation}\label{eq:star}
\|  \Mn(\bar \alpha, \hat \eta)\| \lesssim_{\Pr_n} n^{-1/2}.
\end{equation}
This also implies that
$
\|\Ln(\hat \alpha)\| = \|\Ln(\bar \alpha)\| + o_{\Pr_n}(n^{-1/2}),
$
since  $\|\Ln(\bar \alpha)\| \leq  \|\Ln(\hat \alpha)\|$ and 
$$
\|\Ln(\hat \alpha)\|-  o_{\Pr_n}(n^{-1/2}) \leq \| \Mn(\hat \alpha, \hat \eta) \| \leq   \| \Mn(\bar \alpha, \hat \eta)\| + o(n^{-1/2}) =  \|\Ln(\bar \alpha)\| + o_{\Pr_n}(n^{-1/2}).$$ The former assertion implies that
$
\|\Ln(\hat \alpha)\|^2 = \|\Ln(\bar \alpha)\|^2   + o_{\Pr_n}(n^{-1}),
$
so that
$$
\|\Ln(\hat \alpha)\|^2 - \|\Ln(\bar \alpha)\|^2  =  \|\Gamma_1 (\hat \alpha - \bar \alpha)\|^2 = o_{\Pr_n}(n^{-1}),
$$
from which we can conclude that $\sqrt{n}\| \hat \alpha - \bar \alpha\| \to_{\Pr_n} 0.$

Step 4. (Conclusion).  Given the conclusion of the previous step, the remaining claims are standard and follow
from the Continuous Mapping Theorem and Lemma \ref{CLT:convex}. \qed

\subsection{Proof of Proposition \ref{prop:weighted}}  We have  wp $\to 1$ that, for some constants
$0< u < l < 0$,  $ l \| x\| \leq \| \A x\| \leq u \| x\|$ and  $ l \| x\| \leq \| \hat \A x\| \leq u \| x\|$. Hence \begin{eqnarray*}
&& \sup_{ \alpha \in \mathcal{A} } \frac{\| \hat \A \hat \M^o(\alpha, \hat \eta) - \A \M^o(\alpha, \hat \eta) \| + \|  \A \M^o(\alpha, \hat \eta) - \A\M^o(\alpha,  \eta_0) \|}{r_n + \| \hat \A \hat \M^o (\alpha, \hat \eta)\| + \|\A \M^o(\alpha,  \eta_0)\| }\\
&& \leq  \sup_{ \alpha \in \mathcal{A} } \frac{u}{l} \frac{\|\hat \M^o(\alpha, \hat \eta) - \M^o(\alpha, \hat \eta) \| + \| \M^o(\alpha, \hat \eta) - \M^o(\alpha,  \eta_0) \|}{{(r_n/l) + \|\hat \M^o (\alpha, \hat \eta)\| + \| \M^o(\alpha,  \eta_0)\|  } }  \\
&& \quad +  \sup_{ \alpha \in \mathcal{A} }  \frac{ \| \hat \A - \A\| \|\hat \M^o(\alpha, \hat \eta) \| }{r_n + l \|\hat \M^o (\alpha, \hat \eta)\|  }  \lesssim_{\Pr_n}  o(1) + \| \hat \A - \A\|/l \to_{\Pr_n} 0.
\end{eqnarray*}
The proof that the rest of the conditions hold is analogous and is therefore omitted.\qed

\subsection{Proof of Proposition \ref{prop:one-step}} Step 1.  We define the feasible and infeasible ``one-steps" 
\begin{eqnarray*}
& \check \alpha = \tilde \alpha - \hat F \hat{\M}(\tilde \alpha, \hat \eta), & \hat F = (\hat \Gamma_1' \hat \Gamma_1)^{-1} \hat \Gamma_1',\\
&  \bar \alpha = \alpha_0 - F  \hat{\M}( \alpha_0, \eta_0), & F = (\Gamma_1' \Gamma_1)^{-1} \Gamma_1'.
\end{eqnarray*}
We deduce by (\ref{eq:consistent gamma}) and (\ref{eq:stability})  that
$$
\|\hat F\|\lesssim_{\Pr_n} 1, \quad  \| \hat F \Gamma_1 - I \| \lesssim_{\Pr_n} r_n, \quad  \| \hat F - F\| \lesssim_{\Pr_n} r_n.
$$

Step 2.  By Step 1 and by condition (\ref{eq:derivatives2}), we have that \begin{eqnarray*}
 && \mathsf{D} = \| \hat F \hat{\M}(\tilde \alpha, \hat \eta) - \hat F \hat{\M}(\alpha_0, \eta_0) - \hat F \Gamma_1 (\tilde \alpha - \alpha_0)\| \\
 & &\leq \| \hat F\| \| \hat{\M}(\tilde \alpha, \hat \eta) - \hat{\M}(\alpha_0, \eta_0) - \Gamma_1(\tilde \alpha - \alpha_0) \|\\
&& \lesssim_{\Pr_n} \| \hat{\M}(\tilde \alpha, \hat \eta) - \M (\tilde \alpha, \hat \eta) - \hat{\M}(\alpha_0, \eta_0) \| +
\mathsf{D}_1  \lesssim_{\Pr_n} o(n^{-1/2}) + \mathsf{D}_1,
\end{eqnarray*}
where $\mathsf{D}_1 :=  \| \M (\tilde \alpha, \hat \eta) - \Gamma_1(\tilde \alpha - \alpha_0) \|$.

Moreover,  $\mathsf{D}_1 \leq IV_1 + IV_2 + IV_3$, where wp $\to 1$ by condition (\ref{eq:derivatives2}) and $r_n^2 = o(n^{-1/2})$
\begin{eqnarray*}
&& IV_1: = \| \M(\tilde \alpha, \eta_0) - \Gamma_1(\tilde \alpha - \alpha_0) \| \lesssim \| \tilde \alpha - \alpha_0\|^2 \lesssim r_n^2 = o(n^{-1/2}),\\
&& IV_2:= \|\M(\tilde \alpha, \hat \eta) - \M(\tilde \alpha, \eta_0) - \partial_{\eta'} \M(\tilde \alpha, \eta_0)[\hat \eta - \eta_0] \| \lesssim o(n^{-1/2}),\\
&& IV_3:= \| \partial_{\eta'} \M(\tilde \alpha, \eta_0) [ \hat \eta - \eta_0] \| \lesssim  o(n^{-1/2}).
\end{eqnarray*}
Conclude that $n^{1/2} \mathsf{D} \to_{\Pr_n} 0$.

Step 3.  We have by the triangle inequality and Steps 1 and 2 that
\begin{eqnarray*}
\sqrt{n}\| \check \alpha - \bar \alpha\| && \leq \sqrt{n} \| (I - \hat F \Gamma_1)(\tilde \alpha - \alpha_0)\|
+ \sqrt{n} \| (\hat F - F) \hat{\M}(\alpha_0, \eta_0) \|  + \sqrt{n} \mathsf{D}\\
&&\leq \sqrt{n} \| (I - \hat F \Gamma_1) \| \| \tilde \alpha - \alpha_0\| + \|\hat F - F\| \| \sqrt{n} \hat{\M}(\alpha_0, \eta_0)\| +  \sqrt{n} \mathsf{D}\\
&&\lesssim_{\Pr_n} \sqrt{n}  r_n^2 + o(1) = o(1).
\end{eqnarray*}
Thus, $\sqrt{n}\| \check \alpha - \bar \alpha\|\to_{\Pr_n} 0$, and $\sqrt{n}\| \check \alpha - \hat \alpha\|\to_{\Pr_n} 0$ follows from the triangle inequality and the fact that $\sqrt{n} \| \hat \alpha - \bar \alpha\| \to_{\Pr_n} 0$. \qed

\subsection{Proof of Lemma 2}

The conditions of Proposition 1 are clearly satisfied, and thus the conclusions of Proposition 1 immediately follow.  We also have that, for $\hat \Gamma_1 =\hat \Gamma_1(\hat \eta)$,
$$
\sqrt{n}(\hat \alpha - \alpha_0) = - \hat F \sqrt{n} \hat \M(\alpha_0, \hat \eta),  \quad \hat F =  (\hat \Gamma_1'\hat \Gamma_1)^{-1}  \hat \Gamma_1,
$$
$$
\sqrt{n}(\bar \alpha - \alpha_0) :=  - F \sqrt{n} \hat \M(\alpha_0, \eta_0),  \quad  F =  (\Gamma_1'\Gamma_1)^{-1}  \Gamma_1.$$
We deduce by (33) and (11)  that
$
\|\hat F\|\lesssim_{\Pr_n} 1$ and $\| \hat F - F\| \to_{\Pr_n} 0.$
Hence we have by triangle and H\"{o}lder inequalities and condition (33) that
\begin{eqnarray*}
\sqrt{n}\| \hat \alpha - \bar \alpha\| \leq \| \hat F\| \sqrt{n} \| \hat \M(\alpha_0, \hat \eta)-  \hat \M(\alpha_0, \eta_0))\|
+  \| \hat F - F\| \sqrt{n} \|\hat{\M}(\alpha_0, \eta_0) \| \to_{\Pr_n} 0.
\end{eqnarray*}
The conclusions regarding the uniform validity of inference using $\hat \alpha$, of the form stated
in conclusions of Proposition 2,  follow
from the conclusions regarding the uniform validity of inference using $\bar \alpha$, which
follow from the Continuous Mapping Theorem,  Lemma \ref{CLT:convex}, and the assumed
stability conditions (11). This establishes
the second claim of the Lemma.  Verification of the conditions of Proposition 2
is omitted. \qed

\subsection{Proof of Lemma 3 and 4}The proof of Lemma 3 is given in the main text. As in the
proof of Lemma  3, we can expand: \begin{equation}\label{eq:decompose2}
\sqrt{n}(\hat{\M}_j(\alpha_0, \hat \eta) - \hat{\M}_j(\alpha_0, \eta_0)) =  T_{1,j} + T_{2,j} +T_{3,j},
\end{equation}
where the terms $(T_{l,j})_{l=1}^3$ are as defined in the main text.  We can
further bound $T_{3,j}$ as follows: \begin{equation}\label{eq:term3}
 T_{3,j} \leq T^m_{3,j} + T_{4,j}, \quad \begin{array}{lll}
& &  T^m_{3,j}:=  \sqrt{n} |(\hat \eta - \eta^m_0)' \partial_\eta  \partial_{\eta'}\hat{\M}_j(\alpha_0) (\hat \eta - \eta^m_0)|,\\
& &  T_{4,j}:=  \sqrt{n}  |{\eta^r_0}' \partial_\eta  \partial_{\eta'}\hat{\M}_j(\alpha_0)\eta^r_0|.
\end{array}
\end{equation}
Then $T_{1,j}=  0$ by orthogonality, $T_{2,j} \to_{\Pr_n}0$ as in the proof of Lemma 3.
Since $s^2 \log (pn)^2/n \to 0$, $T^m_{3,j}$ vanishes in probability because, by H\"{o}lder's inequality
and for sufficiently large $n$,
$$
T^m_{3,j} \leq \sqrt{n} \bar T_{3,j} \|\hat \eta - \eta^m_0\|^2 \lesssim_{\Pr_n} \sqrt{n} s \log (pn)/n \to_{\Pr_n} 0.
$$
Also, if $s^2 \log (pn)^2/n \to 0$, $T_{4,j}$ vanishes in probability because, by H\"{o}lder's inequality
and (43),
$$
T_{4,j} \leq \sqrt{n} \|\partial_\eta  \partial_{\eta'}\hat{\M}_j(\alpha_0)\|_{\sf{pw}(\eta_0^r)} \|\eta^r_0\|^2 \lesssim_{\Pr_n} \sqrt{n} s \log (pn)/n \to_{\Pr_n} 0.
$$
The conclusion follows from (\ref{eq:decompose2}). \qed

\subsection{Proof of Lemma 5.} For $m=1,...,k$ and $l=1,...,d$, we can bound each element $\hat{\Gamma}_{1, ml}(\eta)$ of matrix $\hat{\G}_{1}(\eta)$
as follows:
\begin{equation*}
|\hat{\Gamma}_{1, ml}(\hat \eta) - \hat{\G}_{1,ml}(\eta_0)| \leq  \sum_{k=1}^4 T_{k,ml}, \begin{array}{lll}
& &  T_{1,ml}:=  |\partial_{\eta} {\G}_{1,ml}(\eta_0)  '(\hat \eta - \eta_0)|, \\
&& T_{2,ml}:=  |(\partial_\eta \hat{\G}_{1,ml}(\eta_0) -  \partial_\eta  {\G}_{1,ml}(\eta_0) ) '(\hat \eta - \eta_0)|, \\
& &  T_{3,ml}:=   | (\hat \eta - \eta^m_0)' \partial_\eta  \partial_{\eta'}\hat{\G}_{1,ml} (\hat \eta - \eta^m_0)|,\\
& &  T_{4,ml}:= | {\eta^r_0}' \partial_\eta  \partial_{\eta'}\hat{\G}_{1,ml} \eta^r_0|.
\end{array}
\end{equation*}
Under conditions (44) and (45) we have that wp $\to 1$
\begin{eqnarray*}
&& T_{1,ml} \leq \|\partial_{\eta} {\G}_{1,ml}(\eta_0) \|_{\infty} \|\hat \eta - \eta_0\|_{1} \lesssim_{\Pr_n} \sqrt{s^2 \log(pn)/n} \to 0,\\
&& T_{2,ml} \leq \|\partial_\eta \hat{\G}_{1,ml}(\eta_0) -  \partial_\eta  {\G}_{1,ml}(\eta_0) \|_{\infty} \|\hat \eta - \eta_0\|_{1}
\lesssim_{\Pr_n} \sqrt{s^2 \log (pn)/n} \to 0,\\
&& T_{3,ml} \leq   \|\partial_\eta  \partial_{\eta'}\hat{\G}_{1,ml}\|_{\mathsf{sp}(\ell_n s)}\|\hat \eta - \eta^m_0\|^2 \lesssim
_{\Pr_n} s \log (pn)/n \to 0,\\
&& T_{4,ml} \leq  \|\partial_\eta  \partial_{\eta'}\hat{\G}_{1,ml}\|_{\mathsf{pw}(\eta^r_0)}\|\eta^r_0\|^2 \lesssim
_{\Pr_n} s \log (pn)/n \to 0.
\end{eqnarray*}
The claim follows from the assumed growth conditions, since $d$ and $k$ are bounded. \qed

\section{Key Tools}

Let $\Phi$ and $\Phi^{-1}$ denote the distribution and quantile function of $\mathcal{N}(0,1)$.
Note that in particular $\Phi^{-1}(1-a) \leq \sqrt{2 \log (a)}$ for all $a \in (0,1/2)$.

\begin{lemma}[Moderate Deviation Inequality for Maximum of a Vector]\label{Lemma:SNMD} Suppose that
$ \mathcal{S}_{j} :=  {\sum_{i=1}^n U_{ij}}/{\sqrt{ \sum_{i=1}^n U^2_{ij}}},$
where $U_{ij}$ are independent random variables across $i$ with mean zero and finite third-order moments.  Then
$$
\Pr \left( \max_{1 \leq j\leq p }|\mathcal{S}_{j}|  >  \Phi^{-1}(1- \gamma/2p)  \right) \leq \gamma \left(1 + \frac{A}{\ell^3_n}\right),
$$
where $A$ is an absolute constant, provided for $\ell_n > 0$
$$
0 \leq \Phi^{-1}(1- \gamma/(2p))  \leq \frac{n^{1/6}}{\ell_n} \min_{1\leq j \leq p} M^2_{j}-1, \ \  \ M_{j} := \frac{\left( \frac{1}{n} \sum_{i=1}^n \Ep [U_{ij}^2]\right)^{1/2}}{\left(\frac{1}{n} \sum_{i=1}^n \Ep[|U_{ij}|^3] \right)^{1/3}}.
$$
\end{lemma}This result is essentially due to \cite{jing:etal}.  The proof of this result, given in \cite{BCCH12}, follows from a simple combination  of union bounds with their result.  

\begin{lemma}[Laws of Large Numbers for Large Matrices in Sparse Norms]\label{lemma:lln}
Let $s_n$, $p_n$, $k_n$ and $\ell_n$ be sequences of positive constants such that $\ell_n \to \infty$ but $\ell_n/\log n \to 0$
and $c_1$ and $c_2$ be fixed positive constants. Let $(x_i)_{i=1}^n$ be i.i.d. vectors such  that $\| \Ep[x_i x_i'] \|_{\mathsf{sp}(s_n \log n)}\leq c_1$, and
either one of the following holds: (a) $x_i$ is a sub-Gaussian random vector with $\sup_{\|u\| \leq 1 } \| x_i'u\|_{\psi_2,\Pr} \leq c_2$,
where $\|\cdot\|_{\psi_2,\Pr}$ denotes the $\psi_2$-Orlizs norm of a random variable, and $s_n(\log n )(\log (p_n\vee n))/n   \to 0$; or (b) $\| x_i\|_\infty \leq k_n$ a.s. and $k_n^2s_n(\log^4 n) \log(p_n\vee n)/n  \to 0$. Then there is $o(1)$ term such 
that with probability $1- o(1)$:
$$\| \En[x_i x_i'] - \Ep[x_i x_i']\|_{\mathsf{sp}(s_n \ell_n)} \leq o(1), \quad   \| \En[x_i x_i'] \|_{\mathsf{sp}(s_n \ell_n)} \leq c_1 + o(1).$$
\end{lemma}
Under (a) the result follows from Theorem 3.2 in \citen{RudelsonZhou2011} and under (b) the result
follows from \citen{rudelson:vershynin}, as shown in the Supplemental Material of \cite{BC-PostLASSO}.

\begin{lemma}[Useful implications of CLT in $\mathbb{R}^m$]\label{CLT:convex}  Consider a sequence of
random vectors $Z_n$ in $\mathbb{R}^m$ such that $Z_n \leadsto Z = \mathcal{N}(0,I_m)$. The elements of
the sequence and the limit variable need not be defined on the same probability space. Then
$$
\lim_{n \to \infty} \sup_{R \in \mathcal{R}} | \Pb(Z_n \in R) - \Pb(Z \in R) | =0,
$$  
where $\mathcal{R}$ is the collection of all convex sets in $\mathbb{R}^m$.
\end{lemma}

Proof. Let $R$ denote a generic convex set in $\mathbb{R}^m$. Let $R^{\epsilon}=\{ z \in \mathbb{R}^m:  d(z, R) \leq \epsilon\}$
and $R^{-\epsilon}=\{ z \in R:  B(z, \epsilon) \subset R\}$, where $d$ is the Euclidean distance and $B(z, \epsilon) =
\{ y\in \mathbb{R}^m: d(y,z) \leq \epsilon\}$.  The 
set $R^{\epsilon}$ may be empty. By Theorem 11.3.3 in \citen{dudley:rap}, $\epsilon_n: = \rho(Z_n, Z) \to 0$,
where $\rho$ is the Prohorov metric.  The definition of the metric implies that $\Pb(Z_n \in R) \leq \Pb(Z \in R^{\epsilon_n})+ \epsilon_n$. By the
reverse isoperimetric inequality [Prop 2.5. \citen{chenfang}]  $ |\Pb(Z \in R^{\epsilon_n})  -
 \Pb(Z \in R^{}) | \leq m^{1/2} \epsilon_n$. Hence  $\Pb(Z_n \in R) \leq \Pb(Z \in R) + \epsilon_n (1+ m^{1/2})$. Furthermore,
 for any convex set $R$,  $(R^{-\epsilon_n})^{\epsilon_n} \subset R$ (interpreting
 the expansion of an empty set as an empty set). Hence
 for any convex $R$ we have $\Pb(Z \in R^{-\epsilon_n}) \leq \Pb(Z_n \in R) + \epsilon_n$ by definition 
  of Prohorov's metric.  By the reverse isoperimetric inequality $|\Pb(Z \in R^{-\epsilon_n}) - \Pb(Z \in R)| \leq m^{1/2} \epsilon_n$. Conclude    that  $\Pb(Z_n \in R) \geq \Pb(Z \in R) - \epsilon_n (1+ m^{1/2})$. \qed

\footnotesize
\bibliographystyle{aea}
\bibliography{mybib}

\normalsize

\newpage

\bigskip
{\huge Supplemental Appendix for ``Valid Post-Selection and Post-Regularization Inference: An Elementary, General Approach''}\\

In this supplement, we provide proofs for the results from  Section 5 of ``Valid Post-Selection and Post-Regularization Inference: An Elementary, General Approach''.  Equation numbers (1)-(61) refer to equations defined in the main text, and equation numbers (62) and greater are defined in this supplement.

\setcounter{equation}{61}
\setcounter{lemma}{5}

\section{Proof of Proposition 5}

We present the proofs for the case of $p^d =1$; the general case follows similarly.

We proceed to verify the assumptions of Lemma 4 and 5, from which the desired
result follows from Propositions 1 and 2 and Lemma 2.  In what follows, we consider
an arbitrary sequence $\{ \Pr_n \}$ in $\{ \mathbf{P}_n\}$.

Step 1. (Performance bounds for $\hat \eta$).  We noted that Condition AS.1 implies the decomposition (50).
A modification of the proofs of \cite{BCCH12} yields the following
performance bounds for estimator $\hat \eta$ of $\eta_0$:  wp $\to 1$,
\begin{equation}\label{eq:sparse2 in the proof}
\begin{array}{c}
 \|\hat \eta\|_0 \lesssim s,
 \| \hat \eta- \eta^m_0\|_2 \lesssim \sqrt{(s/n) \log (pn)}, \quad   \|\hat \eta - \eta^m_0\|_1
\lesssim  \sqrt{(s^2/n) \log (pn)}.
\end{array}
\end{equation}

Note that the required modification addresses two differences between the development in the present paper and that in \cite{BCCH12}.  First, we impose only that errors are uncorrelated with control regressors and instruments whereas mean independence between errors and controls and instruments is assumed in \cite{BCCH12}. Second, the third step of Algorithm 1 presented in main text requires regressing an estimated response variable on the control regressors.  This extension is handled by noting that the estimation error in the response variable can be treated as additional approximation error in the proofs given in \cite{BCCH12}.  We omit these details for brevity and as they are straightforward.


Step 2. (Preparation). It is convenient to lift the nuisance parameter $\eta$ into a higher dimension and redefine
the signs of its components as follows:
$$
\eta := (\eta_1', \eta_2', \eta_3', \eta'_4, \eta'_5)' := [-\theta', - \vartheta', \gamma', \delta', - \vartheta']'.
$$
With this re-definition, we have
$$
\psi(w_i, \alpha, \eta) = \{ (y_i + x_i' \eta_1) - (d_i + x_i'\eta_2) \alpha \} \{ x_i'\eta_3 + z_i'\eta_4 + x_i'\eta_5\}.
$$
Note also that
$$
\M(\alpha, \eta) =  \Gamma_1(\eta) \alpha + \G_2(\eta), \quad \hat{\M}(\alpha, \eta) = \hat \Gamma_1(\eta) \alpha + \hat \Gamma_2(\eta),
$$ $$
\Gamma_1(\eta) = \Ep[ \partial_{\alpha} \psi(w_i, \alpha, \eta )], \quad \hat \Gamma_1(\eta) = \En[\partial_{\alpha} \psi(w_i, \alpha, \eta ) ].
$$

We  compute the following partial derivatives:
\begin{eqnarray*}
&&\partial_{\eta'} \psi(w_i) := \partial_{\eta'} \psi(w_i, \alpha_0, \eta_0) =  [x_i' \varrho_i, -\alpha_0x_i' \varrho_i, x_i'\varepsilon_i, z_i'\varepsilon_i,  x_i' \varepsilon_i],\\
&& \partial_{\alpha} \psi(w_i, \alpha, \eta ) =  -\{d_i + x_i'\eta_2 \} \{ x_i'\eta_3 + z_i'\eta_4 + x_i'\eta_5\}, \\
&& \partial_{\alpha} \psi(w_i):= \partial_{\alpha} \psi(w_i, \alpha_0, \eta_0) = -\rho^d_i \varrho_i,\\
&&
\partial_{\eta'}\partial_{\alpha} \psi(w_i):= \partial_{\eta'}\partial_{\alpha} \psi(w_i, \alpha_0, \eta_0) =  [0, -x_i' \varrho_i, -x_i' \rho^d_i, -z_i' \rho^d_i, -x_i' \rho^d_i]',\\
\end{eqnarray*} \begin{eqnarray*}
&& \partial_{\eta} \partial_{\eta'}  \psi(w_i, \alpha, \eta )  =
\left[
\begin{array}{ccccc}
 0 & 0  & x_ix_i' & x_iz_i' & x_ix_i' \\
0  &  0 &  -\alpha x_ix_i' & -\alpha x_i z_i' & -\alpha x_ix_i' \\
 x_ix_i' &   -\alpha x_ix_i' &   0    &  0  & 0 \\
  z_ix_i'  &  -\alpha z_i x_i'   & 0 & 0 & 0 \\
 x_ix_i' &  -\alpha x_ix_i' & 0 & 0  & 0
\end{array}
\right],
\\
&& \partial_{\eta} \partial_{\eta'} \partial_\alpha \psi(w_i, \alpha, \eta ) =
\left[ \ \ \begin{array}{ccccc}
  0 & 0 & 0 & 0 & 0 \\
 0 &  0 &  -x_i x_i' & -x_i z_i' & -x_ix_i' \\
 0   &  -x_i x_i'  &   0    &  0  & 0 \\
 0 &   -z_i x_i'  & 0 & 0 & 0 \\
 0 & -x_ix_i' & 0 & 0 & 0
\end{array} \ \
\right].
\end{eqnarray*}

Step 3. (Verification of Conditions of Lemma 4).

Application of Lemma 6, condition $\| \alpha_0 \| \leq \sC$ holding by Condition SM, and Condition RF, yields that wp $\to 1$,
\begin{eqnarray*}
&& \sqrt{n}\| \partial_\eta  \M(\alpha_0, \eta_0) - \partial_{\eta} \hat \M(\alpha_0, \eta_0) \|_{\infty} = \|  \sqrt{n} ( \En -\Ep) \partial_\eta \psi(w_i) \|_\infty  \lesssim\\
&  & \quad  \lesssim  \max_{j}  \sqrt{[\En [ (\partial_\eta \psi(w_i))^2_j ]} \sqrt{ \log (pn) }  \lesssim_{\Pr_n}   \sqrt{ \log (pn) }.
\end{eqnarray*}

Application of the triangle inequality, condition $\| \alpha_0 \| \leq \sC$, and Condition RF yields
$$
\| \partial_\eta \partial_{\eta'} \hat \M(\alpha_0, \eta_0) \|_{\mathsf{sp}(\ell_n s) }\leq C  \| \partial_\eta \partial_{\eta'} \En [f_i f_i] \|_{\mathsf{sp}(\ell_n s) } \lesssim_{\Pr_n} 1,
$$
where $C$ depends on $\sC$.

Moreover, application of the triangle inequality and the Markov inequality yields, for any deterministic $a \neq 0$,
$$
\| \partial_\eta \partial_{\eta'} \hat \M(\alpha_0, \eta_0) \|_{\mathsf{pw}(a) } 
\leq C  \| \partial_\eta \partial_{\eta'} \En [f_i f_i] \|_{\mathsf{pw}(a) } \lesssim_{\Pr_n} 1,
$$
where $C$ depends on $\sC$.

We  have by Condition SM and the law of iterated expectations:
\begin{eqnarray*}
&& \Omega = \Ep[ \psi^2(w_i, \alpha_0, \eta_0) ] = \Ep[\varepsilon_i^2 \varrho_i^2] \in \Ep [\varrho_i^2] \cdot [\sc, \sC] \in [\sc^2, \sC^2],
\\
&&  \Ep[ \psi^{\sq/2}(w_i, \alpha_0, \eta_0) ] \leq \Ep[|\varepsilon_i \varrho_i|^{\sq/2}] \leq \sqrt{ \Ep [|\varepsilon_i|^\sq] } \sqrt{\Ep [|\varrho_i|^\sq ]} \leq \sC.
\end{eqnarray*}
 Application of Lyapunov's Central Limit Theorem yields,
 $$
\Omega^{-1/2}    \hat \M(\alpha_0, \eta_0)  \leadsto \mathcal{N}(0, 1).
 $$

Next, $ \hat \Omega(\alpha_0) = \En[ \psi^2(w_i, \alpha_0, \hat \eta) ]
 $ is consistent for $\Omega$.  The proof of this result
follows similarly to the (slightly more difficult) proof of consistency of $ \hat \Omega= \En[ \psi^2(w_i, \hat \alpha, \hat \eta) ]
$ for $\Omega$, which is given below.

 All conditions of Lemma 4 are now verified.

 Step 4.  (Verification of Conditions of Lemma 5).

Application of Lemma 6 and Conditions RF yields that with probability $1- o(1)$,
\begin{eqnarray*}
&& \sqrt{n}\|  \partial_\eta \hat \G_1(\alpha_0, \eta_0) -  \partial_\eta \G_1(\alpha_0, \eta_0) \|_{\infty} = \|  \sqrt{n} ( \En -\Ep) \partial_\eta \partial_\alpha \psi(w_i) \|_\infty  \\
&  & \quad  \lesssim  \max_{j}  \sqrt{[\En [ (\partial_\eta \partial_\alpha \psi(w_i))^2_j ]} \sqrt{ \log (pn) }  \lesssim_{\Pr_n}   \sqrt{ \log (pn) }.
\end{eqnarray*}

Application of the triangle inequalities and Condition RF yields:
$$
\| \partial_\eta \partial_{\eta'} \hat \G_1(\alpha_0, \eta_0) \|_{\mathsf{sp}(\ell_n s) }\lesssim  \| \partial_\eta \partial_{\eta'} \En [f_i f_i \|_{\mathsf{sp}(\ell_n s) } \lesssim_{\Pr_n} 1.
$$

Moreover, application of the triangle inequalities and the Markov inequality yields, for any deterministic $a \neq 0$,
$$
\| \partial_\eta \partial_{\eta'} \hat \G_1(\alpha_0, \eta_0) \|_{\mathsf{pw}(a) } 
\lesssim  \| \partial_\eta \partial_{\eta'} \En [f_i f_i \|_{\mathsf{pw}(a) } \lesssim_{\Pr_n} 1.
$$

Next, by Condition SM we have
$$
\Gamma_1 = \Ep[ \rho^d_i \varrho_i] = \Ep [\varrho^2_i] \in [\sc, \sC].
$$

By Condition RF we have
\begin{eqnarray*}
\|\partial_\eta \G_1(\alpha_0, \eta_0) \|_\infty = \|\Ep[\partial_{\eta'}\partial_{\alpha} \psi(w_i)]\|_\infty & \leq&  \max_{j} \left( \Ep[|f_{ij} \varrho_i|]  \vee \Ep[|f_{ij} \rho^d_i|]  \right) \\
& \leq&  \max_{j}  \left( \sqrt{[\Ep[|f_{ij} \varrho_i|^2]} \vee \sqrt{|\Ep[|f_{ij} \rho^d_i|^2]} \right)  \leq \sqrt{\sC}.
\end{eqnarray*}
This, as well as previous steps verify conditions of the Lemma 5, which are sufficient
to establish that $|\hat \alpha - \alpha_0| \lesssim_{\Pr_n} n^{-1/2}$, which is needed in the last step below.

Next, we show $\hat V_n - V_n \to_{\Pr_n} 0$.  Given the stability conditions
established above, this follows from $\hat \Gamma_1(\hat \eta) - \Gamma_1 \to_{\Pr_n} 0$,
which follows from Lemma 5, and from 
$\hat \Omega -  \Omega \to_{\Pr_n} 0 $. Recall that  $ \hat \Omega= \En[ \psi^2(w_i, \hat \alpha, \hat \eta) ]
$ and let $\hat \Omega_0 = \En[ \psi^2(w_i, \alpha_0, \eta_0) ]$.   Since $\hat \Omega_0 - \Omega \to_{\Pr_n} 0$ by the Markov
inequality, it suffices to show that 
$\hat \Omega -  \hat \Omega_0 \to_{\Pr_n} 0 $. Since $\hat \Omega - \hat\Omega_0 =  (\sqrt{\hat \Omega} - \sqrt{\hat \Omega_0}) 
(\sqrt{\hat \Omega} + \sqrt{\hat \Omega_0})$, it suffices to show that $(\sqrt{\hat \Omega} - \sqrt{\hat \Omega_0}) \to_{\Pr_n} 0$.
 By the triangle inequality and some simple calculations, we have
$$
|\sqrt{\hat \Omega} - \sqrt{\hat \Omega_0}| \lesssim  \mathsf{D} := I_2 I_\infty \sqrt{\En[\varrho_i^4]}
+ I_2 I_\infty II_2 II_\infty + II_2 II_\infty \sqrt{\En[\varepsilon_i^4]},
$$
where the terms are defined below. Let
\begin{eqnarray*} && \hat \varepsilon_i = \hat \rho^y_i - \hat \rho^d_i \hat \alpha,   \ \ \varepsilon_i = \rho^y_i -  \rho^d_i \alpha,\\
&& \hat \varrho_i = x_i'\hat \gamma + z_i'\hat \delta - x_i'\hat \vartheta, \ \ \varrho_i = x_i' \gamma_0 + z_i'\delta_0 - x_i'\vartheta_0.
\end{eqnarray*}
Then
\begin{eqnarray*}
&& |\hat \varepsilon_i - \varepsilon_i | \leq |x_i'(\theta_0- \hat \theta)| + |\rho^d_i (\alpha_0-\hat \alpha)| + |x_i'(\hat \vartheta  - \vartheta_0) \alpha_0| + |x_i'(\hat \vartheta  - \vartheta_0)  (\hat \alpha-\alpha_0)|,\\
&& |\hat \varrho_i - \varrho_i| \leq |z_i'(\hat \delta-\delta_0)| + |x_i'(\hat \gamma - \gamma_0)| + |x_i'(\hat \vartheta- \vartheta_0)|.
\end{eqnarray*}
Then  the terms $I_2, I_\infty, II_2, II_\infty$ are defined and bounded, using elementary inequalities and Condition RF, as follows:
  \begin{eqnarray*}
&& I_2 := \sqrt{ \En[(\hat \varepsilon_i - \varepsilon_i)^2]} \lesssim_{\Pr_n}   \sqrt{ s \log(pn)/n } + \sqrt{ \En[ \rho_i^{d2}]} |\hat \alpha - \alpha_0| \\
&&  \quad  \quad + \sqrt{ s \log(pn)/n } |\alpha_0|+ \sqrt{ s \log(pn)/n } |\hat \alpha - \alpha_0|  \to_{\Pr_n} 0, \\
 && I_\infty: = \max_{i \leq n} |\hat \varepsilon_i - \varepsilon_i|  \lesssim  \max_{ij} |f_{ij}| \sqrt{s^2 \log(pn)/n }
  + \max_{i \leq n} | \rho_i^d|| \hat \alpha - \alpha_0| \\
 & & \quad  \quad +   \max_{ij} |f_{ij}| \sqrt{s^2 \log(pn)/n } | \alpha_0|+ \max_{ij} |f_{ij}| \sqrt{s^2 \log(pn)/n } | \hat \alpha - \alpha_0|  \to_{\Pr_n} 0,\\
&&  II_2 := \sqrt{ \En[(\hat \varrho_i - \varrho_i)^2]} \lesssim_{\Pr_n}   \sqrt{ s \log(pn)/n } \to 0,
 \\
&& II_\infty: = \max_{i \leq n} |\hat \varrho_i - \varrho_i| \lesssim \max_{ij} |f_{ij}| \sqrt{s^2 \log(pn)/n }\to_{\Pr_n} 0,
 \end{eqnarray*}
where we have used the relations,
  $|\hat \alpha - \alpha_0| \lesssim_{\Pr_n} n^{-1/2}$, and
$\En[|\rho^d_i|^{2}]\lesssim_{\Pr_n} 1$,
$ \max_{i \leq n} \rho_i^d \lesssim_{\Pr_n} n^{1/\sq},
$ for $\sq>4$, holding
by Condition SM, and we have used the fact that
$$
\En[(x_i'\{\hat \theta - \theta_0\})^2 +
(x_i'\{\hat \vartheta - \vartheta_0\})^2 +
(x_i'\{\hat \gamma - \gamma_0\})^2+
(z_i'\{\hat \delta - \delta_0\})^2] \lesssim_{\Pr_n}  \frac{s \log(pn)}{n}.
$$
The latter follows from the following argument, for example, 
wp $\to 1$,
\begin{eqnarray*}
&& \En[(x_i'\{\hat \theta - \theta_0\})^2]  \leq 2 \En[(x_i'\{\hat \theta - \theta^m_0\})^2] + 
2 \En[(x_i'\theta^r_0)^2 ]\\
&& \lesssim \| \En[f_i f_i']  \|_{\mathsf{sp}(\ell_n s)}  \| \hat \theta - \theta^m_0\|^2 + 
\| \En[f_i f_i']  \|_{\mathsf{pw}(\theta^r_0)}  \| \theta^r_0\|^2 \\
&&  \lesssim_{\Pr_n} s \log (pn)/n,
\end{eqnarray*}
since  wp $\to 1$  $\|\hat \theta - \theta^m_0\|_0 \leq \ell_n s$,  $ \| \hat \theta - \theta^m_0\|^2 
\lesssim s \log (pn)/n $, $ \| \theta^r_0\|^2 \lesssim s/n$, by Step 1 and Condition AS.1 (see
decomposition (50)), and $ \| \En[f_i f_i']  \|_{\mathsf{sp}(\ell_n s)}  \lesssim_{\Pr_n} 1$ holding
by Condition RF and $\| \En[f_i f_i']  \|_{\mathsf{pw}(\theta^r_0)}  \lesssim_{\Pr_n} 1$
holding by Markov inequality and Condition RF.

Since $\En[\varrho_i^4] + \En[\varepsilon_i^4] \lesssim_{\Pr_n} 1$ by Condition SM, we conclude that $\mathsf{D} \to_{\Pr_n} 0 $.\qed

\footnotesize

\bibliographystyle{aea}
\bibliography{mybib}

\end{document}